\documentclass[11pt]{article}

\usepackage{amsmath, amsfonts, amssymb,latexsym,wasysym,graphicx}
\input epsf.sty

\newtheorem{Theorem}{Theorem}[section]
\newtheorem{Lemma}{Lemma}[section]
\newtheorem{Proposition}[Lemma]{Proposition}

\newtheorem{Definition}[Lemma]{Definition}

\newcommand{\BEQ}{\begin{equation}}     
\newcommand{\BEA}{\begin{eqnarray}}
\newcommand{\BD}{\begin{displaymath}}
\newcommand{\EEQ}{\end{equation}}       
\newcommand{\EEA}{\end{eqnarray}}
\newcommand{\ED}{\end{displaymath}}

\newcommand{\del}{\delta}
\newcommand{\Del}{\Delta}
\newcommand{\eps}{\varepsilon}          

\newcommand{\g}{{\mathfrak{g}}}

\newcommand{\R}{\mathbb{R}}

\newcommand{\N}{\mathbb{N}}

\renewcommand{\P}{\mathbb{P}}

\newcommand{\Id}{{\mathrm{Id}}}

\def\proba{{\mathbb{P}}}
\def\esper{{\mathbb{E}}}
\def\T{{\mathbb{T}}}

\def\Var{{\mathrm{Var}}}

\newcommand{\eop}{\hfill $\Box$}        

\newcommand{\II}{{\rm i}}               
\newcommand{\half}{{1\over 2}}          
\renewcommand{\vec}[1]{\boldsymbol{#1}} 

\catcode`\@=11
\def\numberbysection{\@addtoreset{equation}{section}
        \def\theequation{\thesection.\arabic{equation}}}
\numberbysection

\begin{document}

\vspace*{1.5cm}
\begin{center}
{\Large \bf From constructive field theory to fractional stochastic calculus. (I) An introduction: rough path theory and perturbative heuristics.}

\end{center}

\vspace{2mm}
\begin{center}
{\bf  Jacques Magnen and J\'er\'emie Unterberger}
\end{center}

\vspace{2mm}
\begin{quote}

\renewcommand{\baselinestretch}{1.0}
\footnotesize
{Let $B=(B_1(t),\ldots,B_d(t))$ be a $d$-dimensional fractional Brownian motion
with Hurst index $\alpha\le 1/4$, or more generally a Gaussian process whose paths have the same local regularity. Defining properly iterated integrals of $B$ is a difficult
task because of the low H\"older regularity index of its paths. Yet rough path theory shows it is
the key to the construction of a stochastic calculus with respect to $B$, or to solving differential
equations driven by $B$.

We intend to show in a forthcoming series of papers how to desingularize iterated integrals by a weak singular non-Gaussian perturbation of the Gaussian measure defined by a limit in law
procedure.
 Convergence is
proved by using  "standard" tools of constructive field theory, in particular cluster expansions and renormalization. These powerful tools allow optimal estimates of the moments and
call for an extension of the Gaussian tools such as for instance the Malliavin calculus.

This first paper aims to be both a presentation of the basics of rough path
theory to physicists, and of perturbative field theory to probabilists; it is
only heuristic, in particular because the desingularization of iterated
integrals is really a {\em non-perturbative} effect. It is also meant to be a
general motivating introduction to the subject, with some insights into
quantum field theory and stochastic calculus. The interested reader should
read in a second time the companion article \cite{MagUnt2} 
(or a preliminary version \cite{MagUnt3}) for the
constructive proofs.
}
\end{quote}

\vspace{4mm}
\noindent
{\bf Keywords:}
fractional Brownian motion, stochastic integrals, rough paths, constructive field theory, Feynman diagrams,
renormalization, cluster expansion.

\smallskip
\noindent
{\bf Mathematics Subject Classification (2000):}  60F05, 60G15, 60G18, 60H05, 81T08, 81T18.

\tableofcontents



\section{Introduction}


A major achievement of the probabilistic school since the middle of the 20th century is the
study of diffusion equations, in connection with Brownian motion or more generally
Markov processes -- and also with partial differential equations, through the Feynman-Kac
formula -- with many applications in physics and chemistry \cite{VK}. One of the main tools
is stochastic calculus with respect to semi-martingales $M$. An adapted integral such as $\int_s^t X(u)
dM(u)$ may be understood as a limit in some sense to be defined. Classically one uses
piecewise linear interpolations, $\sum_{s\le t_1<\ldots<t_N\le t} X(t_i)(M(t_{i+1})-
M(t_i))$ or $\sum_{s\le t_1<\ldots<t_N\le t} \frac{X(t_i)+X(t_{i+1})}{2} (M(t_{i+1})-
M(t_i))$; these approximations define in the limit $N\to\infty$
the It\^o, resp. Stratonovich integral. The latter one is actually obtained e.g. 
if $M=W$ is Brownian motion and $X(t)=f(W_t)$ with $f$ smooth   as the limit
$\lim_{\eps\to 0}\int_s^t f(W_{\eps}(u))dW_{\eps}(u)$ for any smooth approximation
$(W_{\eps})_{\eps>0}$ of $W$ converging a.s. to $W$  (see \cite{WZ}, or
\cite{KS} p. 169). The Stratonovich integral
$\int_s^t X(u)d^{Strato} M(u)$ has an advantage over the It\^o integral in that it agrees
with the fundamental theorem of calculus, namely, $F(M(t))=F(M(s))+\int_s^t F'(M(u)) d^{Strato}M(u)$.

\medskip

The semi-martingale approach fails altogether when considering stochastic processes with
lower regularity. Brownian motion, and more generally semi-martingales (up to time reparametrization), are $(1/2)^-$-H\"older,
i.e. $\alpha$-H\"older for any $\alpha<1/2$ \footnote{Recall that a continuous path $X:[0,T]\to\R$ is $\alpha$-H\"older, $\alpha\in(0,1)$, if $\sup_{s,t\in[0,T]} \frac{|X_t-X_s|}{|t-s|^{\alpha}}<\infty$.}.
 Processes with $\alpha$-H\"older paths, where
$\alpha\ll 1/2$, are maybe less common in nature but still deserve interest.
 Among these, the family of multifractional Gaussian processes is perhaps the most
widely studied \cite{PLV}, but one may also cite diffusions on fractals \cite{HL}, sub- or superdiffusions
in porous media \cite{Goldenfeld,Lesne} and the
fascinating multi-fractal random measures/walks in connection with turbulence and
two-dimensional Liouville quantum gravity \cite{BM,DS}. Many models in hydrodynamics take as input
a space-time noise which is often chosen colored in space \cite{Iva}. In this respect, let us
mention in particular the Kraichnan model for passive advection of scalars, for which anomalous
correlation exponents \cite{KupMur,FalGawVer,BerGawKup} may be expanded in $\alpha$ for $\alpha\to 0$.

\medskip

We concentrate in this article on {\em multiscale Gaussian processes}
 (the terminology is ours) with {\em scaling dimension} or
more or less equivalently {\em H\"older regularity} $\alpha\in(0,1/2)$, the best-known example of
which being {\em fractional Brownian motion} (fBm for short) with Hurst index $\alpha$, $B^{\alpha}(t)$ or simply $B(t)$ \footnote{It is (up to a constant) the unique self-similar Gaussian process with stationary increments. The last property implies that its derivative is a (distribution-valued)
 {\em stationary field}.}.
We consider more precisely  a {\em two-dimensional fBm}, $B(t)=(B_1(t),B_2(t))$, with independent, identically distributed components \footnote{The
one-dimensional case is very different and much simpler,  and has been treated in \cite{GNRV}.}. The covariance kernel $\esper B_i(s)B_j(t)=\half\del_{i,j}(|s|^{2\alpha}+|t|^{2\alpha}
-|t-s|^{2\alpha})$ is that of an {\em integrated colored noise} in the physical terminology \footnote{In other words (informally at least) $\esper B'_i(s)B'_j(t)\sim -c_{\alpha}|t-s|^{2\alpha-2}$ instead of
$\del(t-s)$.}. It is a process with {\em long-range, negative correlations}, which is quite
unusual from a statistical physics point of view; but the emphasis here is on the short-distance
(or ultra-violet) behaviour, not on the long-distance one.

 The simplest non-trivial stochastic integral is then
\BEQ {\cal A}(s,t):=\int_s^t dB_1(t_1)\int_s^{t_1} dB_2(t_2)=\int_s^t (B_2(u)-B_2(s))dB_1(u) , \EEQ
a twice iterated integral, where $B=(B_1(t),B_2(t))$ is a two-component fBm with independent,
identically distributed components. Since
$$\int_s^t dB_1(t_1)\int_s^{t_1} dB_2(t_2)+\int_s^t dB_2(t_2)
\int_s^{t_2} dB_1(t_1)=(B_1(t)-B_1(s))(B_2(t)-B_2(s)),$$
 one is mainly interested in
the antisymmetrized quantity (measuring a {\em signed} area, as follows from the Green-Riemann
formula),
\BEA  {\cal LA}(s,t)&:=&\int_s^t dB_1(t_1)\int_s^{t_1} dB_2(t_2)-\int_s^t dB_2(t_2)
\int_s^{t_2} dB_1(t_1) \nonumber\\
&=&\int_s^t (B_2(u)-B_2(s))dB_1(u)- (B_1(u)-B_1(s))dB_2(u), \nonumber\\ \EEA
called {\em L\'evy area}.  The
corresponding Stratonovich integral, obtained as a limit either by linear interpolation or
by more refined Gaussian approximations \cite{CQ02,Nua,Unt08,Unt08b}, has been shown to {\em diverge} as soon as
$\alpha\le 1/4$.

This seemingly no-go theorem, although clear and derived by straightforward computations
that we reproduce in short in section 1, appears to be a puzzle when put in front of the results of {\em rough path theory} \cite{Lyo98,LyoQia02,Gub,Lej,Lej-bis,FV}. The essential
idea conveyed by this theory --  we shall make this precise in section 2 -- is that a path $\Gamma:\R\to\R^d$ with H\"older regularity index $\alpha\in(0,1)$ must be seen as
the {\em projection} onto the $d$ first components of some "essentially arbitrary" {\em rough path over $\Gamma$}, denoted by
 \BEQ {\bf\Gamma}:\R^2\ni (s,t)\mapsto {\bf\Gamma}_{ts}:=({\bf\Gamma}_{ts}^1,\ldots,{\bf\Gamma}_{ts}^N)\in \R^d\times\R^{d^2}\times\ldots\times
\R^{d^N}, \EEQ
 $N=\lfloor 1/\alpha\rfloor$ \footnote{where $\lfloor \ .\ \rfloor$ stands for the integer part of its argument.}, which may be interpreted as iterated integrals of $\Gamma$
in a limiting sense, namely, $\lim_{\eps\to 0} \int_s^td\Gamma_{i_1}^{\eps}(t_1)=\Gamma_{i_1}(t)-\Gamma_{i_1}(s)$ for the $d$ first components and
\BEQ \lim_{\eps\to 0} \int_s^t d\Gamma_{i_1}^{\eps}(t_1)\int_s^{t_1} d\Gamma_{i_2}^{\eps}(t_2),
\ldots,\lim_{\eps\to 0} \int_s^t d\Gamma_{i_1}^{\eps}(t_1)\ldots\int_s^{t_{N-1}}
d\Gamma_{i_n}^{\eps}(t_n) \EEQ
for the remaining ones, for some smooth family of approximations $(\Gamma^{\eps})_{\eps>0}$ of $\Gamma$. The limit must be understood in
a  H\"older norm sense, as explained in section 2. In other words, there exist infinitely many different families of approximations of $B$ leading to as many different definitions of its iterated integrals! Alas, Gaussian approximations are unfortunately seemingly unable to produce such a definition for fBm with Hurst index $\alpha\le 1/4$.

\bigskip

Our project in this series of papers is to define an explicit rough path over fBm with arbitrary Hurst index, or
more generally {\em multiscale Gaussian fields} (see  a companion article
\cite{MagUnt2} and \cite{MagUnt3} for a preliminary version   ) by an explicit,
probabilistically meaningful
limiting procedure, thus solving at last the problem of constructing a full-fledged,
Stratonovich-like integration with respect to fBm.

Let us explain our strategy for $1/6<\alpha<1/4$. Roughly speaking, our rough path is obtained by making $B=(B(1),B(2))$ interact through a
weak but singular quartic, non-local interaction, which plays the r\^ole of a squared
kinetic momentum, or {\em bending energy}, and makes its L\'evy area -- and at the same time the iterated integrals of higher order --
 {\em finite}. Following the common use
of quantum field theory, this is implemented by multiplying (probabilists would say: {\em penalizing}) the Gaussian measure by
the exponential weight $e^{-\half c'_{\alpha} \int\int {\cal L}_{int}(\phi_1,\phi_2)(t_1,t_2) |t_1-t_2|^{-4\alpha}
dt_1 dt_2}$ \footnote{The unessential constant $c'_{\alpha}$ is fixed e.g. by demanding that the Fourier transform of
 the kernel $c'_{\alpha}|t_1-t_2|^{-4\alpha}$ is the function $|\xi|^{4\alpha-1}$.} , with
\BEQ {\cal L}_{int}(\phi_1,\phi_2)(t_1,t_2) =\lambda^2 \left\{
(\partial {\cal A}^+)(t_1)  (\partial {\cal A}^+)(t_2)+ (\partial {\cal A}^-)(t_1)(\partial {\cal A}^-)(t_2) \right\},  \EEQ
where:  $\lambda$ (the coupling parameter) is a small, positive constant; $\phi_1,\phi_2$ are the (infra-red divergent) stationary fields associated to $B_1,B_2$, with covariance kernel as in eq. (\ref{eq:cov-phi}),
and similarly,  $ {\cal A}^{\pm}$ are stationary {\em left- and right-turning fields}, built
out of $\phi_1,\phi_2$ and representing the singular part of the L\'evy area (see section 1 for details). As usual in quantum field theory, one considers first the truncated
measure obtained by an "ultra-violet cut-off" and on a finite "volume" (or finite horizon, in the probabilistic terminology) $V=[-T,T]$, i.e. one multiplies
 the Fourier transforms of  the fields  $\phi_1,\phi_2$ 
by some cut-off function with compact support in $[-M^{\rho},M^{\rho}]$ (for some fixed constant $M>1$) and integrates over $V$; see Definition \ref{def:cut-off} for the
precise procedure. Then $\partial {\cal A}^{\pm}$ are replaced by the {\em truncated
 quantities} $(\partial {\cal A}^{\pm})^{\to\rho}$ built out of the {\em truncated fields} $\phi^{\to\rho}$.
The truncated interacting Lagrangian reads
\BEA &&\half c'_{\alpha} \int\int_{V\times V}  |t_1-t_2|^{-4\alpha} {\cal L}_{int}^{\to\rho}(\phi_1,\phi_2)
(t_1,t_2) dt_1 dt_2+\int_V {\cal L}^{\to\rho}_{bdry}  \nonumber\\
&& \qquad :=\half c'_{\alpha}\lambda^2 \int\int_{V\times V} |t_1-t_2|^{-4\alpha}
\left\{
(\partial {\cal A}^+)^{\to\rho}(t_1)  (\partial {\cal A}^+)^{\to\rho}(t_2) \right. \nonumber\\
&& \left. \qquad \qquad + (\partial {\cal A}^-)^{\to\rho}(t_1)
(\partial {\cal A}^-)^{\to\rho}(t_2) \right\} dt_1 dt_2 
+\int_V {\cal L}^{\to\rho}_{bdry}, \nonumber\\  \EEA
where $ {\cal L}_{bdry}^{\to\rho}$ is some singular ``Fourier boundary term'' multiplied by
an evanescent factor $M^{-\kappa\rho}$ ($\kappa>0$), which cures unwanted difficulties due to
the ultra-violet cut-off \footnote{The exact form of ${\cal L}_{bdry}^{\to\rho}$ requires detailed
constructive explanations and will not be required here. It is to be found in the companion
article \cite{MagUnt2}.}.
 When $\rho$ and $V$ are finite, the underlying Gaussian fields are smooth, which ensures the existence of the penalized measure.  The assertion is that {\em the penalized measures
converge weakly when $\rho,|V|\to\infty$ to some well-defined, unique  measure}, while the {\em truncated iterated integrals} themselves
converge in law to a {\em rough path over $B$}.

\medskip

Note that the statistical weight is maximal when 
 $\partial{\cal A}^+=
\partial {\cal A}^-=0$, i.e. for sample paths which are ``essentially''  straight lines. Another way to
motivate this interaction (following an image due to A. Lejay) 
 is to understand that the divergence of the L\'evy area is due to the accumulation in a small region of space of small loops \cite{Lej}; the statistical
weight is unfavorable to such an accumulation. On the other hand, the law of the quantities in the first-order Gaussian chaos, characterized by  the $n$-point functions 
\BEA && \langle
B_{i_1}(x_1)\ldots B_{i_n}(x_n)\rangle_{\lambda} \nonumber\\
&& \quad =\frac{1}{Z}\esper\left[B_{i_1}(x_1)\ldots B_{i_n}(x_n)
e^{-\half c'_{\alpha}  \int\int  {\cal L}_{int}(\phi_1,\phi_2)(t_1,t_2) |t_1-t_2|^{-4\alpha} dt_1 dt_2 
} \right], \nonumber\\ \EEA
$i_1,\ldots,i_n=1,2$, where 
\BEQ  Z:=
\esper\left[
e^{-\half c'_{\alpha}  \int\int  {\cal L}_{int}(\phi_1,\phi_2)(t_1,t_2) |t_1-t_2|^{-4\alpha} dt_1 dt_2 } \right] \EEQ
is a normalization constant playing the r\^ole of a {\em partition function}, 
  is insensitive to the interaction \footnote{In the two preceding equations, $\esper\left[ \ \cdot\ e^{-\half c'_{\alpha} \int\int  {\cal L}_{int}(\phi_1,\phi_2)(t_1,t_2) |t_1-t_2|^{-4\alpha} dt_1 dt_2 }\ \right]$ stands for the limit of $\esper \left[ \ \cdot\ e^{-\half c'_{\alpha} \int\int_{V\times V}  {\cal L}_{int}^{\to\rho}(\phi_1,\phi_2)(t_1,t_2) |t_1-t_2|^{-4\alpha} dt_1 dt_2
  +\int_V {\cal L}^{\to\rho}_{bdry} }\ \right]$ when $\rho,|V|\to\infty$ as we explained above.}. Thus we have built a rough path {\em over fBm}. This conveys the idea that the paths have been straightened by removing
  {\em in average} small bubbles of scale $M^{-\rho}$. In doing so, the paths of the limiting process when $\rho\to\infty$ are indistinguishable from those of $B$, {\em but} higher-order integrals
  have been corrected so as to become finite.

Starting from the above field-theoretic description, the proof of finiteness and H\"older regularity of the L\'evy area for $\lambda>0$ small enough follows, despite some specific features, the broad
scheme of constructive field theory, see e.g. the monographies \cite{Abd,Mas,Riv}.
Constructive field theory is a program originally advocated in the sixties by
A. S. Wightman   \cite{Wigh},  the aim of which was to give explicit
examples of field theories with a non-trivial interaction ; see Glimm and Jaffe's book \cite{GJb} for an introduction and references therein for an extensive bibliography.
 Let us  give  a short guide to the history of the subject.

 The first contribution was made in 1965 by E. Nelson who introduced a scale analysis \cite{Nel}
to control the divergence of a model whose only divergence comes from Wick ordering.
J. Glimm and A. Jaffe  introduced the phase space analysis \cite{GJa} for  models having a finite number of divergent graphs. 
The cluster expansion was devised by  J. Glimm, A. Jaffe and T. Spencer \cite{GJS} to control infinite volume limits.

The Roman team  \cite{Benfa} realized that the above phase space analysis was in some sense a continuous space version of  the  block-spin expansion, first written by Kadanoff
for the Ising model, and then  made into a major tool both in high-energy and statistical physics by K. G. Wilson through the introduction of the concept of {\em renormalization group} \cite{Wilb,Wilc}.
The  multiscale expansion was devised  in the eighties  to provide a rigorous version of Wilson's renormalization group e.g. including the flow of the
effective parameters: see 
\cite{GKc} for the block-spin approach, and \cite{FMRS} for the continuous space {\em multi-scale cluster expansion}. 

For some  fermionic theories a simpler version of these constructions is available, due
to the fact that (contrary to the bosonic case)
 the series expansion  in terms of the effective coupling constants is {\em convergent}
  \cite{Les,AbdRiv3}. 
The multi-scale cluster expansion  has also allowed to  study models  with a singularity around a surface,
 like the so-called {\em jellium model} of interacting, non-relativistic fermions
  \cite{FMRTa,DMR}, modelling the generation of Cooper pairs, in connection with the famous
BCS (Bardeen-Cooper-Schrieffer) theory of
   supraconductivity \cite{FMRTb}.

In this work we use the multi-scale cluster expansion developed in \cite{FMRS} more than twenty years ago which seems to us 
 the most appropriate for these probabilistic models; it reduces to the minimum the use of abstract combinatorial identities and algebra, to the benefit of a very intuitive and visual (though sometimes
heavy) tree  expansion.

\bigskip

The main theorem may be stated as follows. As a rule, we denote in this article by $\esper[...]$ the Gaussian expectation and
by $\langle ...\rangle_{\lambda,V,\rho}$
 the expectation with respect to the $\lambda$-weighted interaction measure
  with scale $\rho$ ultraviolet
cut-off restricted to a compact interval $V$,
 so that in particular $\esper[...]=\langle ...\rangle_{0,\infty}$.

\begin{Theorem} \label{th:0.1}

Assume $\alpha\in(\frac{1}{6},\frac{1}{4})$. Consider for $\lambda>0$ small enough the family of probability measures
 (also called: {\em $(\phi,\partial\phi,\sigma)$-model})
\BEQ\proba_{\lambda,V,\rho}(\phi_1,\phi_2)=
e^{-\half c'_{\alpha} \int\int dt_1 dt_2 |t_1-t_2|^{-4\alpha} {\cal L}^{\to\rho}_{int}(\phi_1,\phi_2)(t_1,t_2)-\int {\cal L}_{bdry}^{\to\rho}} d\mu^{\to\rho}(\phi_1)
d\mu^{\to\rho}(\phi_2),\EEQ
where $d\mu^{\to\rho}(\phi_i)=d\mu(\phi_i^{\to\rho})$ is a Gaussian measure obtained  by an ultra-violet cut-off
at Fourier momentum $|\xi|\approx M^{\rho}$ ($M>1$), see Definition \ref{def:cut-off}. Then $(\P_{\lambda,V,\rho})_{V,\rho}$
 converges in law when $|V|,\rho\to\infty$ to some   measure
$\P_{\lambda}$, and the associated  iterated integrals
 $$ \int_s^t d\phi_{i_1}^{\to\rho}(t_1)\int_s^{t_1}d\phi_{i_2}^{\to\rho}(t_2),\ldots,\int_s^t d\phi_{i_1}^{\to\rho}(t_1)\int_s^{t_1} d\phi_{i_2}^{\to\rho}(t_2)\ldots\int_s^{t_{n-1}} d\phi_{i_n}^{\to\rho}(t_n),\ldots$$ 
 converge in law   to a rough path over $B$.
\end{Theorem}

\bigskip

The result is not difficult to understand {\em heuristically},
 at least for quantum field theory experts, if one resorts to  the non-rigorous perturbation
theory (see sections 3 and 4). First, by a  Hubbard-Stratonovich transformation (a functional
Fourier transform), one replaces the non-local
interaction ${\cal L}(\phi_1,\phi_2)(t_1,t_2)|t_1-t_2|^{-4\alpha}$
 with a local interaction ${\cal L}(\phi_1,\phi_2,\sigma)(t)$ depending
on a two-component exchange particle field $\sigma=(\sigma_+(t),\sigma_-(t))$.
Then a {\em Schwinger-Dyson identity} (a functional integration by parts) relates the moments of $\cal A$ to those of $\sigma$. Simple power-counting arguments show that a connected $2n$-point function of $\sigma$ alone is superficially divergent if and only if $1-4n\alpha\ge 0$. Thus, restricting to $\alpha>1/8$, one only needs to renormalize the two-point function. Since the renormalized propagator of $\sigma$ is screened by a positive, infinite mass term,
 the theory is free once one has integrated out the $\sigma$-field, hence one retrieves the underlying Gaussian theory $(\phi_1,\phi_2)$. The Schwinger-Dyson identity then shows that the two-point functions of $\cal A$ have been made finite. Finally, simple arguments (not developed here) yield the convergence of higher-order
iterated integrals in the interacting theory provided $\alpha>1/6$.

\bigskip

Whereas these heuristic arguments are not difficult to follow in principle, they do not constitute at all a proof. Theorem 0.1 is proved in the companion article \cite{MagUnt2} by following --
as explained above -- the general scheme of constructive field theory. Although the constructive method is really a multi-scale refinement of the previous arguments, explaining it precisely
is actually a formidable task, which is in general very much model-dependent, whereas perturbative renormalization always follows more or less the same lines; briefly said, the difference lies
in the difference between a formal power series expansion and an analytic proof of convergence for a given quantity. This task we perform at long length and in great generality in the
companion article, with the view of making constructive arguments into classical  mathematical tools which probabilists may eventually reemploy.    
\bigskip

Here is an outline of the article. 

We begin in Section 1 by recalling classical arguments (due to the second author) explaining the divergence of the L\'evy area for $\alpha\le 1/4$, which is the starting point for all the
story \cite{Unt08}; Fourier normal ordering \cite{Unt-Holder,Unt-fBm} -- an indispensable tool for the sequel -- is introduced there. Section 2 is a brief introduction into rough path theory, mainly for
non-experts. Subsections 2.1 and 2.2 are standard and may be skipped by experts, whereas subsection 2.3 -- a brief summary of the previous contributions of the second author to the subject --
gives the context in which this series of papers arose.

The heart of the article is Section 3 and Section 4. Our problem is recast into a quantum field theoretic language in section 3; we take the opportunity to explain the basis of quantum field
theory and renormalization at the same time. The interaction term is introduced at this point, where it comes out naturally. Finally, section 4 is dedicated to a heuristic perturbative
"proof" of the convergence of the L\'evy area of the interacting process, and serves also in some sense as an introduction to the companion paper \cite{MagUnt2}.

\bigskip


\section{A Fourier analysis of the L\'evy area}


The quantity we want to define in the case of fractional Brownian motion is the following.

\begin{Definition}[L\'evy area]

The L\'evy area of a two-dimensional path $\Gamma:\R\to\R^2$ between $s$ and $t$ is
the area between the straight line connecting $(\Gamma_1(s),\Gamma_2(s))$ to $(\Gamma_1(t),\Gamma_2(t))$ and the curve $\{(\Gamma_1(u),\Gamma_2(u)); s\le u\le t\}$.
It is given by the following antisymmetric quantity,
\BEQ {\cal LA}_{\Gamma}(s,t):=\int_s^t d\Gamma_1(t_1) \int_s^{t_1} d\Gamma_2(t_2) -
\int_s^t d\Gamma_2(t_2) \int_s^{t_2} d\Gamma_1(t_1).\EEQ

\end{Definition}

The purpose of this section is to show by using Fourier analysis why the L\'evy area of fBm 
 diverges when $\alpha\le 1/4$. This is hopefully understandable to physicists, and also profitable to probabilists who
are aware of other proofs of this fact, originally proved in \cite{CQ02}, because Fourier analysis is essential in the analysis of Feynman graphs which shall be needed in section 4. We follow here the computations made in
\cite{Unt-fBm} or \cite{Unt-ren}.

\begin{Definition}[Harmonizable representation of fBm]

Let $W(\xi),\xi\in\R$ be a complex Brownian motion \footnote{Formally,
$\langle W'(\xi_1)W'(\xi_2)\rangle=0$ and
$\langle W'(\xi_1)\overline{W'(\xi_2)}\rangle=\del(\xi_1-\xi_2)$
if $\xi_1,\xi_2>0$.} such that $W(-\xi)=-\overline{W(\xi)}$, and
\BEQ B_t:=(2\pi c_{\alpha})^{-\half} \int_{-\infty}^{+\infty} \frac{e^{\II t\xi}-1}{\II\xi} |\xi|^{\half-\alpha} dW(\xi),
\quad t\in\R.\EEQ
\end{Definition}

The field $B_t, t\in\R$ is called {\em fractional Brownian motion}
\footnote{The constant $c_{\alpha}$ is  conventionally chosen
 so that $\esper (B_t-B_s)^2=|t-s|^{2\alpha}$.} . Its paths are almost surely
$\alpha^-$ H\"older, i.e. $(\alpha-\eps)$-H\"older for every $\eps>0$. It has dependent
but identically distributed (or in other words, stationary) increments $B_t-B_s$.
In order to gain translation invariance, we shall rather use the closely related {\em stationary
process}
\BEQ \phi(t):= \int_{-\infty}^{+\infty} \frac{e^{\II t\xi}}{\II\xi} |\xi|^{\half-\alpha} dW(\xi),
\quad t\in\R\EEQ
-- with covariance
\BEQ \langle \phi(x)\phi(y)\rangle =  \int e^{\II \xi (x-y)} \frac{1}{|\xi|^{1+2\alpha}} d\xi
\label{eq:cov-phi} \EEQ
--
which is  {\em infrared divergent}, i.e. {\em divergent around $\xi=0$}. However,  the increments $\phi(t)-\phi(s)=B_t-B_s$ are well-defined for any $(s,t)\in\R^2$.

In order to understand the analytic properties of the L\'evy area of fBm, we shall resort to
a Fourier transform. One obtains, using the harmonizable representation of fBm,
\BEA  {\cal A}(s,t) &:=& \int_s^t dB_1(t_1) \int_s^{t_1} dB_2(t_2) \nonumber\\
&=& \frac{1}{2\pi c_{\alpha}}\int \frac{dW_1(\xi_1) dW_2(\xi_2)}{|\xi_1|^{\alpha-1/2} |\xi_2|^{\alpha-1/2}} \int_s^t dt_1 \int_s^{t_1} dt_2 \ \cdot\ e^{\II (t_1\xi_1+t_2\xi_2)}.
\EEA

The L\'evy area ${\cal LA}(s,t):={\cal LA}_B(s,t)$ is obtained from this twice iterated integral by
antisymmetrization. Note that ${\cal LA}(s,t)$ is homogeneous of degree $2\alpha$ in $|t-s|$ since $B(ct)-B(cs)$, $c>0$ has same law as $c^{\alpha}(B(t)-B(s))$ by
self-similarity.

Expanding the right-hand side yields an expression which is not homogeneous in $\xi$.
Hence it is preferable to define instead the following stationary quantity called {\it skeleton integral}, which depends only
on {\em one} variable,
\BEA {\cal A}(t) &:=& \int^t dB_1(t_1)\int^{t_1} dB_2(t_2)\nonumber\\
 &=& \frac{1}{2\pi c_{\alpha}}
\int \frac{dW_1(\xi_1) dW_2(\xi_2)}{|\xi_1|^{\alpha-1/2} |\xi_2|^{\alpha-1/2}} \int^t dt_1 \int^{t_1} dt_2 \ \cdot\ e^{\II (t_1\xi_1+t_2\xi_2)} \nonumber\\
&=& \frac{1}{2\pi c_{\alpha}} \int \frac{dW_1(\xi_1) dW_2(\xi_2)}{|\xi_1|^{\alpha-1/2} |\xi_2|^{\alpha-1/2}}
\cdot\ \frac{e^{\II t(\xi_1+\xi_2)}}{[\II(\xi_1+\xi_2)][\II\xi_2]}, \EEA
where by definition $\int^t e^{\II u\xi} du=\frac{e^{\II t\xi}}{\II\xi}.$ From
${\cal A}(t)$ and the one-dimensional skeleton integral
\BEQ \phi_i(t)=(2\pi c_{\alpha})^{-\half} \int^t dB_i(u)=\int \frac{dW_i(\xi)}{|\xi|^{\alpha-1/2}} \ \cdot\ \frac{e^{\II t\xi}}{\II \xi},\EEQ
which is the above-defined infra-red divergent stationary process associated to $B$, 
one easily retrieves ${\cal A}(s,t)$ since
\BEA {\cal A}(s,t) &=& \int_s^t dB_1(t_1) \left( \int^{t_1} dB_2(t_2) -\int^s dB_2(t_2)\right) \nonumber\\
&=& {\cal A}(t)-{\cal A}(s)+{\cal A}_{\partial}(s,t),  \label{eq:dec1} \EEA
where $(2\pi c_{\alpha})^{\half} {\cal A}_{\partial}(s,t):=(B_1(t)-B_1(s))\phi_2(s)$ (called {\em boundary term})
is a {\em product of first-order integrals}.

One may easily estimate these quantities in each sector $|\xi_1|\gtrless|\xi_2|$. In
practice, it turns out that estimates are easiest to get {\em after} a permutation of the
integrals (applying Fubini's theorem) such that (for twice or multiple iterated integrals
equally well) {\em innermost (or rightmost) integrals bear highest Fourier frequencies};
this is the essence of {\em Fourier normal ordering} \cite{Unt-Holder,FoiUnt,Unt-resume}. This gives a somewhat different
decomposition with respect to (\ref{eq:dec1}) since $\int_s^t dB_1(t_1)\int_s^{t_1}dB_2(t_2)$ is rewritten as $-\int_s^t dB_2(t_2)\int_t^{t_2} dB_1(t_1)$
in the "negative" sector $|\xi_1|>|\xi_2|$. After some elementary computations, one gets
the following.

\begin{Lemma}
Let
\BEA {\cal A}^+(t)&:=&2\pi c_{\alpha} \int^t dt_1 \int^{t_1} dt_2 {\cal F}^{-1}\left(
(\xi_1,\xi_2)\mapsto {\bf 1}_{|\xi_1|<|\xi_2|} ({\cal F}B'_1)(\xi_1) ({\cal F}B'_2)(\xi_2)
\right) (t_1,t_2) \nonumber\\
&=&  \int_{|\xi_1|<|\xi_2|} \frac{dW_1(\xi_1)dW_2(\xi_2)}{|\xi_1|^{\alpha-1/2} |\xi_2|^{\alpha-1/2}} \ \cdot\ \frac{e^{\II t(\xi_1+\xi_2)}}{[\II(\xi_1+\xi_2)][\II\xi_2]}
\EEA
and
\BEA {\cal A}^-(t)&:=& 2\pi c_{\alpha} \int^t dt_2 \int^{t_2} dt_1 {\cal F}^{-1}\left(
(\xi_1,\xi_2)\mapsto {\bf 1}_{|\xi_2|<|\xi_1|} ({\cal F}B'_1)(\xi_1) ({\cal F}B'_2)(\xi_2)
\right) (t_1,t_2) \nonumber\\
&=& \int_{|\xi_2|<|\xi_1|} \frac{dW_1(\xi_1)dW_2(\xi_2)}{|\xi_1|^{\alpha-1/2} |\xi_2|^{\alpha-1/2}} \ \cdot\ \frac{e^{\II t(\xi_1+\xi_2)}}{[\II(\xi_1+\xi_2)][\II\xi_1]}.
\EEA

Then
\BEQ {\cal A}(s,t)=\frac{1}{2\pi c_{\alpha}} \left\{ ({\cal A}^+(t)-{\cal A}^+(s))-({\cal A}^-(t)-{\cal A}^-(s))+
({\cal A}^+_{\partial}(s,t)-{\cal A}^-_{\partial}(s,t)) \right\}, \label{eq:1.8} \EEQ
the {\em boundary term} ${\cal A}^+_{\partial}-{\cal A}^-_{\partial}$ being given by
\BEA && {\cal A}^+_{\partial}(s,t)-{\cal A}^-_{\partial}(s,t)= \left\{ -
\int_{|\xi_1|<|\xi_2|} \frac{(e^{\II t\xi_1}-e^{\II s\xi_1})e^{\II s\xi_2}}{[\II \xi_1][\II \xi_2]} \right. \nonumber\\
&&\left. \qquad \qquad +\int_{|\xi_2|<|\xi_1|} \frac{(e^{\II t\xi_2}-e^{\II s\xi_2})e^{\II t\xi_1}}{[\II \xi_1][\II \xi_2]} \right\}\ \cdot\   \frac{dW_1(\xi_1)dW_2(\xi_2)}{|\xi_1|^{\alpha-1/2} |\xi_2|^{\alpha-1/2}}. \nonumber\\ \EEA

\end{Lemma}

Two lines of computations show immediately that 
\BEA  \Var {\cal A}^{\pm}_{\partial}(s,t) &\lesssim& \int |e^{\II t\xi}-e^{\II s\xi}|^2 |\xi|^{-1-4\alpha}
d\xi \nonumber\\
&\lesssim& \int_{|\xi|>\frac{1}{|t-s|}} \frac{d\xi}{|\xi|^{1+4\alpha}}+\int_{|\xi|<\frac{1}{|t-s|}} \frac{|t-s|^2
|\xi|^2}{|\xi|^{1+4\alpha}} d\xi \nonumber\\
&\lesssim&  |t-s|^{4\alpha}, \label{eq:1.13} \EEA
so that (essentially by the Kolmogorov-Centsov lemma, see section 2) the H\"older regularity indices of $B_1$ and $B_2$
add in the case of the boundary term, to produce a quantity which is $2\alpha^-$-H\"older.
(Note that  the artificial infrared divergence at $\xi_1=0$ disappears when Taylor expanding $e^{\II t\xi_1}-e^{\II s\xi_1}$). On the other hand, letting $\xi:=\xi_1+\xi_2$ and introducing an ultra-violet cut-off at $|\xi_2|=\Lambda\gg 1$, one may see
for instance ${\cal A}^+(t)$ as an inverse random Fourier transform of the integral
$\xi\mapsto\int_{|\xi-\xi_2|<|\xi_2|}^{\Lambda}
 \frac{dW_2(\xi_2)}{\xi_2} \frac{1}{|\xi-\xi_2|^{\alpha-1/2} |\xi_2|^{\alpha-1/2}}$,
whose variance diverges like $\int^{\Lambda} \frac{d\xi_2}{\xi_2^{4\alpha}}=O(\Lambda^{1-4\alpha})$
or $O(\ln\Lambda)$  in the ultra-violet limit $\Lambda\to\infty$
as soon as $\alpha\le 1/4$. Note that the ultraviolet divergence is in the region $|\xi_1|,|\xi_2|\gg
|\xi|$. 

\bigskip

It is apparent that the central r\^ole in this decomposition is played by the Fourier
projection operator $D({\bf 1}_{|\xi_1|<|\xi_2|})={\cal F}^{-1} \left( {\bf 1}_{|\xi_1|<|\xi_2|} \ \cdot\ {\cal F}(\ .\ ) \right).$ Since ${\cal A}_{\partial}^{\pm}$
are obtained by Fourier projecting $(B_1(t)-B_1(s))\phi_2(s)$, or $(B_2(t)-B_2(s))\phi_1(t)$, which are perfectly well-defined products of continuous
fields \footnote{apart from the spurious infra-red divergence (see above)}, it
was clear from the onset that these would be regular terms. Hence singularities come only
from the {\em one-time quantity}  ${\cal A}^{\pm}(t)$, which {\em does not} split into a
 product of first-order integrals, and that we shall call the {\em singular part of the L\'evy area}.


\section{An introduction to rough paths}


\subsection{General issues}


Let $\Gamma=(\Gamma_1(t),\ldots,\Gamma_d(t))$ be a smooth path with $d$ components. As explained in the Introduction,
 the L\'evy area of $\Gamma$, ${\bf\Gamma}^2_{ts}(i,j):=\int_s^t d\Gamma_i(t_1)\int_s^{t_1}
d\Gamma_j(t_2)$ is the simplest non-trivial iterated integral of $\Gamma$. The interest for iterated integrals of $\Gamma$ comes from the study of two closely related problems in the case when $\Gamma$
is not regular any more.

\medskip

1.\ {\bf Integration along an irregular path.}  \medskip

Assume one wants to define the integral of the (say, smooth) one-form $f:=\sum_{j=1}^d f_j(x) dx^j$ along the path $\Gamma$, namely, the quantity\\
 $\int_s^t fd\Gamma:=\sum_j \int_s^t f_j(\Gamma(u))d\Gamma_j(u)$.
Since $\Gamma$ is not differentiable, $d\Gamma_j(u)$ may not be understood as $\frac{d\Gamma_j}{du}\cdot du$, and the very meaning of this quantity is unclear. Unfortunately, the Riemann-type sum
$\sum_j \sum_{i=0}^{n-1} f_j(\Gamma(t_i))(\Gamma_j(t_{i+1})-\Gamma_j(t_i))$, with $s=t_0<\ldots<t_i=s+\frac{i}{n}(t-s)<\ldots<t_n=t$, may be shown to diverge in general as soon as $\alpha\le \half$ \footnote{From
a naive bound by the $1$-variation of the path, $\sum_j\sum_i |\Gamma_j(t_{i+1})-\Gamma_j(t_i)|=O(n^{1-\alpha})$, one would come to the erroneous conclusion that the Riemann-type sums diverge when $\alpha<1$. The so-called {\em Young theory of integration} (see e.g. \cite{Lej}) lowers the barrier to $\alpha=\half$ by taking into account the H\"older regularity of the integrand $f(\Gamma(t))$.}. 

A Taylor expansion to order $N$ of the integrand yields (coming back to the case of a {\em regular path}) the improved Riemann-type sum
\BEQ \sum_{i=0}^{n-1} \sum_{p=1}^N \sum_{j_1,\ldots,j_p=1}^d \frac{\partial^{p-1} f_{j_p}}{\partial x_{j_1}\ldots\partial x_{j_{p-1}}}(\Gamma(t_i)){\bf\Gamma}^p_{t_{i+1},t_i}(j_1,\ldots,j_p), \label{eq:Riemann}\EEQ
where 
\BEQ {\bf\Gamma}^p_{t_{i+1},t_i}(j_1,\ldots,j_p)=\int_s^t d\Gamma_{j_1}(t_1)\ldots\int_s^{t_{p-1}}d\Gamma_{j_p}(t_p) \label{eq:bf-Gamma} \EEQ
 is a $p$-th order iterated integral. The problem is, if $\Gamma$ is {\em irregular}, iterated
integrals of $\Gamma$ are a priori ill-defined for the same reasons as before.

\medskip
2.\ {\bf Solutions of differential equations driven along an irregular path}

Consider the differential equation
\BEQ dy_t=\sum_{i=1}^d V_j(y(t))d\Gamma_j(t). \label{eq:eq-dif} \EEQ
The following series gives a formal solution,
\BEQ  y_t=y_s+\sum_{N=1}^{\infty} \sum_{1\le i_1,\ldots,i_N\le d} [V_{i_1}\cdots
V_{i_N}\cdot {\mathrm{Id}}](Y_s) \ \cdot\ {\bf \Gamma}^{ts}(i_1,\ldots,i_N),
  \label{eq:series} \EEQ
with ${\bf\Gamma}$ as in eq. (\ref{eq:bf-Gamma}). Solutions are usually computed by using some
iterated numerical scheme. For instance, the Euler scheme of rank $N$ gives the solution to
(\ref{eq:eq-dif}) as the limit when $n\to\infty$ of the compound mapping,
\BEQ   
\Phi({\bf X}_{t,t_{n-1}}; \cdots \Phi({\bf X}_{t_2,t_1};\Phi({\bf X}^{t_1,s};y_s)\cdots), \label{eq:Euler} \EEQ
where $\Phi({\bf \Gamma}_{ts};y_s)$ is the series (\ref{eq:series}) truncated to order $N$.
If one takes for $\Gamma$ an $\alpha$-H\"older path, one stumbles again into the same problem of
defining ${\bf\Gamma}_{ts}=({\bf\Gamma}^1_{ts},\ldots,{\bf\Gamma}^N_{ts})$.

\bigskip

In both cases, the hope is that, if one finds some (non necessarily unique!)
 way of defining iterated integrals of $\Gamma$ with the correct regularity properties, then the refined Riemann-type sums (\ref{eq:Riemann}) or
Euler scheme (\ref{eq:Euler}) converge when the mesh $\frac{t-s}{n}$ goes to $0$. Rough path theory shows
this is possible \footnote{Furthermore, the limit is $\alpha$-H\"older and satisfies nice continuity
properties with respect to the path $\Gamma$.} provided one chooses $N\ge \lfloor 1/\alpha\rfloor$ --
 here we choose $N=\lfloor 1/\alpha\rfloor$ minimal -- and 
\BEQ {\bf\Gamma}_{ts}=\left({\bf\Gamma}^1_{ts}(i_1)_{1\le i_1\le d},\ldots,
{\bf\Gamma}^N_{ts}(i_1,\ldots,i_N)_{1\le i_1,\ldots,i_N\le d}\right):\R^2\to \R^d\times\R^{d^2}\times\ldots\times
\R^{d^N} \EEQ is a rough path with H\"older regularity index $\alpha$ in the following sense:

\begin{Definition}[rough path]

 \label{def:rough-path}

An {\em $\alpha$-H\"older continuous  rough path over $\Gamma$} is a functional ${\bf\Gamma}^n_{ts}(i_1,\ldots,i_n)$, $n\le \lfloor N:=
1/\alpha\rfloor$, $i_1,\ldots,i_n\in\{1,\ldots,d\}$, such that ${\bf\Gamma}_{ts}(i)=
\Gamma_t(i)-\Gamma_s(i)$ are the increments of $\Gamma$, and the following 3 properties are satisfied:

\begin{itemize} \item[(i)] (H\"older continuity) ${\bf\Gamma}^n_{ts}(i_1,\ldots,i_n)$ is $n\alpha$-H\"older
continuous as a function of two variables, namely, $\sup_{s,t\in\R}
\frac{|{\bf\Gamma}^n_{ts}(i_1,\ldots,i_n)|}{|t-s|^{\alpha}}<\infty.$

 \item[(ii)] (Chen property) 
\BEA && 
{\bf\Gamma}^n_{ts}(i_1,\ldots,i_n)={\bf\Gamma}^n_{tu}(i_1,\ldots,i_n)+  {\bf\Gamma}^n_{us}(i_1,\ldots,i_n)+ \nonumber\\
&& \qquad 
\sum_{n_1+n_2=n} {\bf\Gamma}^{n_1}_{tu}(i_1,\ldots,i_{n_1}){\bf\Gamma}^{n_2}_{us}(i_{n_1+1},\ldots,i_{n});\EEA

\item[(iii)] (shuffle property) \BEQ
    {\bf\Gamma}_{ts}^{n_1}(i_1,\ldots,i_{n_1}){\bf\Gamma}^{n_2}_{ts}(j_1,\ldots,j_{n_2})=\sum_{\vec{k}\in
    Sh(\vec{i},\vec{j})} {\bf\Gamma}^{n_1+n_2}_{ts}(k_1,\ldots,k_{n_1+n_2}),\EEQ where $Sh(\vec{i},\vec{j})$ --
    the set of shuffles of the words $\vec{i}$ and $\vec{j}$ -- is the subset of permutations of the
    union of the lists $\vec{i},\vec{j}$ leaving unchanged the order of the sublists $\vec{i}$ and
    $\vec{j}$. For instance, ${\bf\Gamma}^2_{ts}(i_1,i_2){\bf\Gamma}^1_{ts}(j_1)=
    {\bf\Gamma}^3_{ts}(i_1,i_2,j_1)+{\bf\Gamma}^3_{ts}(i_1,j_1,i_2)+{\bf\Gamma}^3_{ts}(j_1,i_1,i_2).$
    \end{itemize}
    
    A {\em formal} rough path over $\Gamma$ is a functional satisfying all the above properties
    {\em except} H\"older continuity (i).

\end{Definition}

In a random setting, the H\"older continuity estimates (i) are generally proved as a consequence of moment estimates such as
$\esper |{\bf\Gamma}^n_{ts}|^{2p}\le C_p|t-s|^{2pn\alpha}$, $p\ge 1, n=1,\ldots,N$. This may be seen as a consequence of the well-known
Kolmogorov-Centsov lemma stating that (for a measurable process random $\Gamma$)
{\small \BEQ \left( \esper[|\Gamma(t)-\Gamma(s)|^{2p}]\le C|t-s|^{1+2p\alpha} \right)\Rightarrow \left( \forall\alpha^-<\alpha,
\esper\left[ \left( \sup_{s,t\in[0,T]} \frac{|\Gamma(t)-\Gamma(s)|}{|t-s|^{\alpha^-}} \right)^{2p}\right]<\infty \right) \EEQ }
or more precisely of an extension (or a variant)  of these estimates adaptated to functions of two variables (such as $(s,t)\mapsto {\bf\Gamma}_{ts}^n$)
due to Garsia, Rodemich and Rumsey  \cite{Gar}.

\medskip

In particular, if $\Gamma$ is smooth, then  its natural iterated integrals \\
$\int_s^t d\Gamma_{i_1}(t_1)\ldots\int_s^{t_{n-1}} d\Gamma_{i_n}(t_n)$  satisfy properties (ii) and (iii).

\bigskip

However, it is not clear a priori {\em in what sense} abstract data as in Definition \ref{def:rough-path} should represent iterated integrals in the usual sense.


\subsection{Geometric approach}


The answer to this question comes from a reinterpretation of rough paths in terms of group theory and geometric structures. We generally refer to the book by P. Friz and N. Victoir
\cite{FV} for this paragraph. Consider the signature ${\bf \Gamma}_{ts}=({\bf
\Gamma}_{ts}^1,{\bf \Gamma}^2_{ts},\ldots)$ of a smooth path $\Gamma$  as
\BEQ {\bf\Gamma}(s,t):=1+\sum_{i_1} {\bf\Gamma}^1_{ts}(i_1)X^{i_1}+\sum_{i_1,i_2} {\bf\Gamma}^2_{ts}(i_1,i_2)X^{i_1}\otimes X^{i_2}+\ldots, \EEQ
sitting inside the tensor algebra $T\R^d=\oplus_{n\ge 0}T^n\R^d$, with $X^1,\ldots,X^d$ generating a basis of $\R^d\simeq \T^1\R^d$. Note that the Chen property is trivially equivalent to the property ${\bf\Gamma}(s,t)={\bf\Gamma}(s,u)\otimes
{\bf\Gamma}(u,t)$, implying that ${\bf\Gamma}(s,t)={\bf\Gamma}(0,s)^{\otimes -1}\otimes {\bf\Gamma}(0,t)$ is a {\em multiplicative increment}. In the particular case when $\Gamma(t)=tV$,
$V\in\R^d$ is a straight line, ${\bf\Gamma}(0,t)=\exp t\sum_{i=1}^d V_i X^i$ belongs to $\exp T^1\R^d$. Easy arguments due to Chow show then that $t\mapsto {\bf\Gamma}(0,t)$ is a 
$G$-valued path, where $\g=Lie(G)$ is the {\em free Lie algebra in $d$ generators}, generated as a vector space by the successive commutators $X^{i_1},[X^{i_1},X^{i_2}],[X^{i_1},[X^{i_2},
X^{i_3}]],\ldots$  In rough path theory, one quotients out by $\oplus_{n\ge N+1} T^n \R^d$. Then the quotient Lie algebra $\g_N$ is the {\em free $N$-step nilpotent Lie algebra in
$d$ generators}, and $G_N=\exp \g_N$ is a {\em Carnot group}. When $d=N=2$, $\g_2\simeq \langle X,Y,Z:=[X,Y]\rangle$ is isomorphic to the Heisenberg algebra, and the defect of additivity
of the L\'evy area ${\cal LA}_{\Gamma}(s,t)=\int_s^t d\Gamma_1(t_1)\int_s^{t_1}d\Gamma_2(t_2)-\int_s^t d\Gamma_2(t_2)\int_s^{t_2}d\Gamma_1(t_1)$, measured by the difference
\BEA &&  {\cal LA}_{\Gamma}(s,t)-{\cal LA}_{\Gamma}(s,u)-{\cal LA}_{\Gamma}(u,t)= \nonumber\\
&& \qquad \qquad (\Gamma_1(t)-\Gamma_1(u))(\Gamma_2(u)-\Gamma_2(s))-(\Gamma_2(t)-\Gamma_2(u))(\Gamma_1(u)-\Gamma_1(s)),\nonumber\\ \EEA
is encoded into the non-commutativity of the product in the Heisenberg group, given by (in the exponential coordinates)
 $(x_1,y_1,z_1)\cdot(x_2,y_2,z_2)=(x_1+y_1,x_2+y_2,z_1+z_2+\half(x_1 y_2-x_2 y_1)).$

 Carnot groups are naturally equipped by homogeneous norms coming from the sub-Riemannian Carnot-Carath\'eodory metric induced by horizontal geodesics, i.e.minimizing curves  with tangent vectors in
the Euclidean space  $T^1\R^d$. Then an $\alpha$-H\"older rough path over an $\alpha$-H\"older path $\Gamma$ is simply  an $\alpha$-H\"older $G_N$-valued path (in geometric terms, an $\alpha$-H\"older section of the principal bundle $\R\times G_N$)
 which projects onto $\Gamma$.

One has the  following two fundamental results.
 
 \begin{Proposition} (see  Lyons \cite{Lyo98}, Lyons-Victoir \cite{LyoVic07}, Friz-Victoir \cite{FV})
 
Let $0<\alpha^-<\alpha<1$. 
 
 \begin{enumerate}
 \item (Existence theorem) There exists a (highly non-unique) $\alpha^-$-rough path over any $\alpha$-H\"older path. In geometric
terms, one may lift an $\alpha$-H\"older section of the quotient bundle $\R\times \left(G_N/\exp \oplus_{n=2}^N T^n\R^d\right) \simeq \R\times T^1\R^d
\simeq \R\times \R^d$ into an $\alpha^-$-H\"older section of the principal bundle $\R\times G_N$.

  \item (Approximation theorem) Every $\alpha$-H\"older rough path over $\Gamma$ is the limit in $\alpha^-$-H\"older norm of the corresponding stack of natural iterated integrals over some smooth approximation
  family $\Gamma^{\eps},\eps\to 0$ of $\Gamma$.
 \end{enumerate}
 
 \end{Proposition}

The approximation theorem is essential in that it reduces differential equations driven by $\alpha$-H\"older paths (through a limiting procedure which is often very subtle) to
ordinary differential equations. Estimates for solutions in a deterministic setting are given in full details in the book by
P. Friz and N. Victoir (see \cite{FV}, Chap. 10).

\bigskip
This general approach is however insufficient for many purposes. Drawbacks are of two types:

-- the arguments leading to the existence and approximation theorems are abstract, the first theorem relying on the axiom of choice (due
to the arbitrariness of the lift), and the second one on an interpolation by sub-Riemannian
geodesics which are notoriously complicated objects;

-- in a random setting, this approach produces in principle deterministic, pathwise estimates, which moreover do not depend on the choice of rough path. Even in combination with probabilistic
tools such as the Malliavin calculus, despite beautiful achivements in the case $\alpha>1/4$ 
(such as global existence of solutions for bounded potentials \cite{FV}, existence of a density \cite{CassFri}, ergodicity \cite{Hai}),...) generalizing results known
in the case of diffusion equations,  it does not permit -- in the case of stochastic differential equations driven by  fBm for instance -- to produce anything really better than a local existence
theorem for solutions beyond the barrier $\alpha=1/4$.

\bigskip

Let us mention briefly  {\em en passant} another related approach due to M. Gubinelli \cite{Gub} and called {\em algebraic integration}. Without being too precise, it states the existence of a class of {\em $\Gamma$-controlled paths} -- stable under functional transformations and under  integration along $\Gamma$, and to which solutions of differential
equations driven by $\Gamma$ belong -- whose increments are of the form
\BEQ z_t-z_s=\sum_{n=1}^N \sum_{i_1,\ldots,i_n} \zeta^n_s(i_1,\ldots,i_n) {\bf\Gamma}^n_{ts}(i_1,\ldots,i_n) \label{eq:Gamma-control} \EEQ
for some functions $\zeta^n(i_1,\ldots,i_n)$, up to a remainder $\rho_{ts}$ such that $\rho_{ts}=O(|t-s|^{1+\eps}$, with
$\eps>0$. The right-hand side of (\ref{eq:Gamma-control}) -- viewed as a function
of $t$ -- is a linear combination of the components of the rough path $\Gamma$, while the remainder is sufficiently regular so that conventional estimates apply. This essentially
avoids the use of smooth approximations and requires only the knowledge of the quantities ${\bf\Gamma}^n_{ts}$, $n\le N$.


\subsection{Fourier normal ordering}

In contrast with this geometric approach, the point of view developed by  the second  author is that a rough path over an irregular path $\Gamma$ is something "essentially arbitrary", 
and that one should rather look for {\em explicitly} constructed rough paths with "good" properties, which allow better estimates than the general ones.

Let us summarize very roughly the results obtained so far in the following Proposition:

\begin{Proposition} (see \cite{Unt-Holder,FoiUnt,Unt-ren}) \label{prop:Unt}

\begin{enumerate}
\item A rough path is uniquely determined by an algorithm called {\em Fourier normal ordering algorithm} from its {\em tree data}, which are generalized Fourier normal ordered skeleton
integrals on domains indexed by trees. As a consequence, any {\em arbitrary set of tree data} produces a {\em formal rough path} (see Definition \ref{def:rough-path}).

\item Tree data yielding {\em H\"older-continuous rough paths} by Fourier normal ordering may be obtained by various, explicit regularization schemes applied to Fourier normal ordered tree
skeleton integrals, using multi-scale methods and inspired by the renormalization of Feynman graphs. In particular, one may construct rough paths ${\bf B}=({\bf B}^1_{ts},
\ldots,{\bf B}^N_{ts})$ over fBm such that ${\bf B}^j_{ts}$ is in the $j$-th chaos of fBm \footnote{i.e. may be written as a $j$-linear integral expression in terms of $B$.}.

\end{enumerate}

\end{Proposition}

Fourier normal ordering consists as in section 2 in (1) cutting iterated integrals
like $I^{ts}_{\Gamma}(1,\ldots,n):=\int_s^t d\Gamma_1(t_1)\int_s^{t_1} d\Gamma_2(t_2)\ldots \int_s^{t_{n-1}} d\Gamma_n(t_n)$
into $n!$ pieces by applying the Fourier projection operators ${\cal P}^{\sigma}:=D({\bf 1}_{|\xi_{\sigma(1)}|<\ldots<|\xi_{\sigma(n)}|})$, where $\sigma$ ranges in the group of permutations of $\{1,\ldots,n\}$; (2) rewriting each piece ${\cal P}^{\sigma}I^{ts}_{\Gamma}(1,\ldots,n)$
as a Fourier normal ordered integral over the inverse image of the simplex $\{t>t_1>\ldots
>t_n>s\}$ by $\sigma$ by using Fubini's theorem. The inverse image of the simplex decomposes
as a union of elementary domains indexed by trees \footnote{Given a rooted tree  with $n$ vertices
indexed by $1,\ldots,n$,
one integrates over the domain with coordinates $t_1,\ldots,t_n\in[s,t]$ such that $t_i<t_j$
whenever the vertex $i$ is above the vertex $j$. When the tree is simply a trunk tree with
no branching, one gets a usual iterated integral of order $n$.}.

\smallskip 

Thus the r\^ole of Fourier normal ordering is twofold: (1) it allows a general algebraic (combinatorial) classification of (formal) rough paths; (2) it induces a correct addition of
the H\"older regularity indices of the tree data when recombining them by the Fourier normal ordering algorithm. We have seen an example of this when we estimated the variance of the
boundary terms ${\cal A}_{\partial}^{\pm}$ in section 1.

\medskip

The rough paths described in the above Proposition, in the case of fBm, say, are not obtained by an explicit limiting procedure; yet they suggest very strongly that the construction
of rough paths is closely related to renormalization in quantum field theory. The purpose of the present series of articles is to give a probabilistic construction coming {\em directly}
from quantum field theory. We actually conjecture that (some of) the rough paths of the above Proposition may be obtained by some limiting procedure from the construction
of the next sections.

\section{Definition of the interaction}

We recall that $\int_s^t dB_1(t_1)\int_s^{t_1}dB_2(t_2)$ represents the
area between the straight line connecting $(B_1(s),B_2(s))$ to $(B_1(t),B_2(t))$ and the curve. If the curve turns right, resp. left, then the
L\'evy area increases, resp. decreases. We have seen that  ${\cal A}^{\pm}$ represents in some
sense the {\em singular part} of the L\'evy area.  

\medskip

It is conceivable that $B_1,B_2$ or $\phi_1,\phi_2$ represent the idealized, strongly
self-correlated motion in $\R^2$ of a particle, which -- although rotation-invariant --
may not (probably as a consequence of a mechanical or electromagnetic rigidity
due to the macroscopic dimension of the particle, or any other similar phenomenon)
turn absolutely freely. A natural quantum field theoretic description of this rigidity
phenomenon is to add an interaction Lagrangian of the form ${\cal L}_{int}=(\partial
{\cal A}^{\pm})^2$.  The fundamental intuition here is that the field $B$ is in some sense a {\em mesoscopic field},
while ${\cal A}^{\pm}$ depends on {\em microscopic details of the theory}.

\medskip
 This is explained in great accuracy in
\cite{Lej-bis}, in a mathematical language.
 A. Lejay shows how a path $\Gamma$  may be modified by inserting microscopic bubbles all  along,
resulting in the limit in  a path which is indistinguishable from the original one, while the L\'evy area has been corrected
by an {\em arbitrary} amount.  Let us give a very simple example. Take for $\Gamma$ a straight line $\Gamma(t)=\left(\begin{array}{c}
t\\ 0\end{array}\right)$, and insert (somewhat artificially) microscopic bubbles of size $\eps=M^{-\alpha\rho}$ (covered in a time
$O(M^{-\rho})$) at times which are multiples of $M^{-\rho}$. Then the resulting path $\Gamma^{\eps}$ has a L\'evy area of order
$M^{\rho}\cdot (M^{-\alpha\rho})^2\to_{\rho\to\infty}\infty$, while $\Gamma^{\eps}\to\Gamma$ in $\alpha^-$-H\"older norm
whenever $\alpha^-<\alpha$ since $\frac{|(\Gamma^{\eps}(t)-\Gamma^{\eps}(s))-(\Gamma(t)-\Gamma(s))|}{|t-s|^{\alpha^-}}=O(M^{-(\alpha
-\alpha^-)\rho})\to_{\rho\to\infty} 0$. The inverse process of {\em removing} microscopic bubbles of a given path so as to make 
its L\'evy area finite is of course much more hazardous, and looks a little bit like an ``inverse Joule expansion''  (i.e. like
putting back all the molecules of a gas into the left compartment of a container after removing the wall which separated it from
the right compartment, a statistical physicist's nightmare, sometimes called "Maxwell's devil").

\bigskip
Summarizing the above discussion,  one must search for an interaction which cures the ultra-violet divergences
of the microscopic scale,  without modifying the theory at mesoscopic scale.
This is where quantum field theory comes into play. The interested reader may refer to several excellent treatises on the subject (see e.g. \cite{PesSch} or \cite{LeB}). It is impossible to
give here a self-contained introduction to this theory which is one of the main foundations of the modern physics of both high-energy particles and condensed matter. Let us however explain in an informal way the most essential concepts, and introduce some useful terminology, in order to  fill in the gap between probability theory and physics. We have tried to make the  next two definitions
as precise and as general as possible. In our case the space-time dimension $D$ is simply one.

\begin{Definition}[ultra-violet cut-off] \label{def:cut-off}

\begin{enumerate}
\item Let $M>1$ be a constant, and $\chi^0:\R^D\to\R$, resp. $\chi^1$ a non-negative, compactly supported function such that $\chi^0\equiv 1$ in a neighbourhood of $0$, resp. $\chi^1\equiv 0$
in a neighbourhood of $0$ and $\chi^1\equiv 1$ in a neighbourhood of the hypersquare $\sup_{j=1,\ldots,D}|\xi_j|=1$. These two functions may be chosen such that $(\chi^0,(\chi^j)_{j\ge 1})$, with $\chi^j:=\chi^1(M^{-j}\cdot)$, define a partition of unity, i.e. $\chi^0+\sum_{j\ge 1}\chi^j\equiv 1$.  Let $\rho\in\N$. Then the {\em ultra-violet
cut-off at scale $\rho$} of a function $f:\R^D\to\R^d$ is $f^{\to\rho}:={\cal F}^{-1}\left(\xi\mapsto \left[\sum_{j=0}^{\rho} \chi^j(\xi)\right] {\cal F}f(\xi)\right)$, where $\cal F$ is the Fourier transformation.
Roughly speaking, the ultra-violet cut-off cuts away Fourier components of {\em momentum} $\xi$ such that $|\xi|>M^{\rho}$.
\item Let $C_{\phi}(x,y):=C_{\phi}(x-y)$ be the covariance of a stationary Gaussian field $\phi:\R^D\to\R$. Then $\phi$ has same law as the series  of independent Gaussian fields
$\sum_{j\ge 0} \phi^j$, where $\phi^j$ has covariance kernel
$C_{\phi}^j:={\cal F}^{-1}\left(\xi\mapsto \chi^j(\xi){\cal F}C_{\phi}(\xi)\right)$. The {\em ultra-violet cut-off at scale $\rho$} of the Gaussian field $\phi$ is then $\phi^{\to\rho}:=\sum_{j=0}^{\to\rho} \phi^j$,
with covariance $C_{\phi}^{\to \rho}:=\sum_{j=0}^{\rho} C_{\phi}^j$.
\end{enumerate}

\end{Definition}

Note that (at least for a good choice of the functions $\chi^0,\chi^1$) the Fourier transform of $\phi^j$
is supported on the union of two dyadic slices, $M^{j-1}<|\xi|<M^{j+1}$.
In principle one may extend this decomposition to {\em negative scale indices} $j$, so that the limit $j\to -\infty$ describes the correlations at large distances. In our model however -- and this makes it very different with respect to classical models in statistical physics, see comments below -- it is only the transition from the microscopic scale $\rho$ to the mesoscopic scale which is
non-trivial, and one may essentially restrict to positive indices $j$.

\begin{Definition}[interacting fields]

Let $\phi:\R^D\to\R^d$ be a vector-valued Gaussian process on $\R^D$, $\vec{\lambda}:=(\lambda_1,\ldots,
\lambda_q)$ a set of real parameters, 
 and $P_1,\ldots,P_q$ ($q\ge 1$) homogeneous polynomials on $\R^d\times (\R^d)^D$. Then the {\em interacting theory
with interaction Lagrangian} ${\cal L}_{int}(\phi)(x)=\sum_{p=1}^q \lambda_p P_p(\phi(x);\nabla \phi(x))$ is (provided it exists!) the weak limit $\proba(d\phi)$ of the {\em penalized measures}
\BEQ \proba_{\vec{\lambda}, V, \rho}(d\phi):=\frac{1}{Z_{V,\rho}} e^{-\int_V {\cal L}_{int}(\phi^{\to\rho})(x)dx} d\mu^{\to\rho}(\phi\big|_V),\EEQ
when the volume $|V|$ and the ultra-violet scale $\rho$ go to infinity, where: $V\subset\R^D$ is compact; $d\mu^{\to\rho}(\phi\big|_V)$ is the Gaussian measure corresponding to the cut-off
field $\phi^{\to\rho}$ restricted to the finite volume $V$; $Z_{V,\rho}$ is a normalization constant called {\em partition function} by reference to Gibbs measures.

\end{Definition}

In general, $\phi$ is stationary, which accounts for the {\em finite volume cut-off $V$}, and $\int_V {\cal L}_{int}(\phi^{\to\rho})(x)dx$ diverges when $\rho\to\infty$, which accounts for
the {\em ultra-violet cut-off at scale $\rho$}. The parameters $\lambda_1,\ldots,\lambda_q$ are called {\em bare coupling constants}. Usually the inverse of the covariance kernel of $\phi$ is a differential
operator of the form $C_{\phi}^{-1}=\lambda_{\nabla}\nabla^2+m^2$, where $m$ is called the {\em mass}. (In the case of our model, $C_{\phi}^{-1}$ contains a {\em fractional} derivative operator
instead, but the present
discussion remains valid). Formally (forgetting about the cut-offs) $d\mu(\phi)$ gives the trajectories a weight proportional to the Onsager-Machlup functional $e^{-\half( (\lambda_{\nabla}\nabla^2+m^2)\phi,\phi)}$,
so the parameters $\lambda_{\nabla}$ and $m^2$ play a r\^ole similar to the coupling constants $\lambda_1,\ldots,\lambda_q$, and the sum of the {\em interaction Lagrangian} and of the {\em Onsager-Machlup 
functional}
is called simply the {\em Lagrangian}.

In general also, $\phi$ is self-similar (or at least asymptotically self-similar at short distances), so the term in the Lagrangian $P_p(\phi(x),\nabla\phi(x))dx$ has a certain {\em degree
of homogeneity} with respect to a change of scale $x\mapsto ax$ or equivalently $\xi\mapsto a^{-1}\xi$ after a Fourier transform,
 which gives the main behaviour at large momenta $\xi$ -- or equivalently at short distances --  of the correlations (or so-called $n$-point correlation functions)
$\langle \phi_{i_1}(x_1)\ldots\phi_{i_n}(x_n)\rangle_{V,\rho}:=\int \phi_{i_1}(x_1)\ldots\phi_{i_n}(x_n) \proba_{\vec{\lambda},V,\rho}(d\phi)$.

\bigskip

Here we  take a high-energy physics point of view. Then the {\em bare scale} is $\rho$; in other words, one uses a cut-off at short distances of order $M^{-\rho}\to_{\rho\to\infty} 0$,  and wants to understand the behaviour of the correlations at macroscopic distances \footnote{In statistical physics, the size of the lattice usually gives an explicit cut-off, so one may take $\rho=0$.} . In principle, the theory is hopelessly divergent in the
limit $\rho\to\infty$ if this degree of homogeneity is {\em negative} (the so-called {\em non-renormalizable case}). On the contrary, expanding the exponential $e^{-\int P_p(\phi(x),\nabla\phi(x))
dx}$ into a series leads to only a finite number of diverging terms (called diverging {\em Feynman diagrams})  if the degree of homogeneity
is {\em positive} (the so-called {\em super-renormalizable case}). When this degree of homogeneity is zero (the so-called {\em just renormalizable case}, often the most interesting one in 
practice) closer inspection is needed. In all cases, for a large variety of models, one obtains by  iterated integration with respect to highest Fourier scales (i.e. with respect to the field components
$\phi^{\rho},\phi^{\rho-1},\ldots,\phi^{j+1}$) an effective theory at scale $j$ which may be described in terms of the same Lagrangian  {\em but} with so-called
{\em renormalized parameters}, by opposition to the {\em bare parameters},
$\lambda_p\rightsquigarrow \lambda_p^j$ or $\lambda_{\Del}\rightsquigarrow \lambda_{\nabla}^j$, $m^2\rightsquigarrow (m^2)^j$. One obtains in general a flow for the parameters, i.e. equations of the type $(\lambda_{\nabla}^j,(m^2)^j;(\lambda_p^j), p'=1,\ldots,q):=F(\lambda_{\nabla}^{j+1},(m^2)^{j+1};(\lambda_{p'}^{j+1}),p'=1,\ldots,q)$. Solving this flow
down to small values of $j$ is then the main task of renormalization. An interesting case is when one may show that the contribution of the renormalized vertex
$\lambda_p P_p(\phi,\nabla\phi)$ goes to zero at    distances which are large with respect to the bare scale; then the theory is said to be {\em asymptotically free at large distances}. The best-known examples of this behaviour are maybe
the weakly avoiding path  or the $\phi^4$-theory, both in dimension $D=4$ ;
see  \cite{MI,GKc,FMRS} for rigorous results \footnote{On the other hand, in high-energy physics, the main example in this respect
is that of asymptotic freedom at {\em short distances} (or equivalently {\em high energy})  of quarks \cite{GW,P}, so exactly the opposite point of view with respect to the one we adopt here.}.
Our model is original for it combines in some sense features of models of  both high-energy physics and  statistical physics: namely, the bare scale is $O(M^{-\rho})$, {\em but}
the theory is asymptotically free at large distances. Letting $\rho\to\infty$, the interaction disappears {\em at all finite scales}, hence one retrieves in the end a Gaussian theory, in
which, however, the singular part of the L\'evy area has been cancelled.

\medskip
Perturbative methods are by far the most common in physics, because they are accessible to non-experts. They rely on an asymptotic analysis of the quantities obtained by expanding into
power series in the coupling constants  the exponential weight $e^{-\int {\cal L}_{int}(\phi)(x)dx}$. These are conventionally represented as {\em Feynman graphs} (we shall show some of these later on for our model).
Unfortunately, in all interesting cases, the series diverges by and large because of huge combinatorial factors, hence perturbative theory has only a {\em heuristic status}. Constructive
methods, on the other hand (when they work!), are based on particularly clever {\em finite Taylor expansions}, scale after scale, and produce {\em converging series} (but
not power series!); in other terms, they are
rigorous. However, the technical apparatus needed to explain constructive field theory is much more sophisticated.

\bigskip

Let us now return to the discussion of our model after this long parenthesis. In order to keep track of the degree of homogeneity of the fields -- and to obtain eventually
the expected H\"older regularity indices for iterated integrals -- we need here 
a  {\em just renormalizable}
theory (or, in other terms, an integrated  interaction which is homogeneous of degree $0$).
 Since $(\partial A^{\pm})^2$ is homogeneous of
degree $(4\alpha-2)$ in time, one shall use in fact a {\em non-local} interaction
lagrangian, $\half c'_{\alpha} \int\int |t_1-t_2|^{-4\alpha} {\cal L}_{int}(\phi_1,\phi_2)(t_1,t_2) dt_1 dt_2$, where
\BEQ {\cal L}_{int}(\phi_1,\phi_2)(t_1,t_2)= \lambda^2\left\{ \partial {\cal A}^+(t_1)\partial {\cal A}^+(t_2)+ \partial {\cal A}^-(t_1)\partial {\cal A}^-(t_2)\right\},\EEQ
which is {\em positive} for $\alpha<1/4$ since the kernel $|t_1-t_2|^{-4\alpha}$ is
locally integrable  and positive definite. Thus the Gaussian measure is {\em penalized} by
the singular exponential weight $e^{-\frac{c'_{\alpha}}{2} \int\int {\cal L}_{int}(\phi_1,\phi_2)(t_1,t_2) |t_1-t_2|^{-4\alpha}
dt_1 dt_2}$. Equivalently, using the so-called Hubbard-Stratonovich transformation
\footnote{which is nothing else but an infinite-dimensional extension of the Fourier
transform $\esper[e^{\II\lambda X}]=e^{-\sigma^2\lambda^2/2}$ for  a random variable
$X\sim {\cal N}(0,\sigma^2)$}, we introduce two independent exchange particle fields
$\sigma_{\pm}=\sigma_{\pm}(t)$ with covariance kernel $C_{\sigma_{\pm}}(s,t)=C_{\sigma_{\pm}}(t-s)=\esper \sigma_{\pm}(s)
\sigma_{\pm}(t)=c'_{\alpha} |s-t|^{-4\alpha}$ and rewrite (letting $d\mu(\phi)$, resp.  $d\mu(\sigma)$ be the Gaussian
measure associated to $\phi$, resp.  $\sigma=(\sigma_+,\sigma_-)$)
the partition function $Z:=Z(\lambda)$,
\BEQ Z:=\int e^{-\frac{c'_{\alpha}}{2} \int\int_{\R^2} |t_1-t_2|^{-4\alpha} {\cal L}_{int}(\phi_1,\phi_2)(t_1,t_2)
dt_1 dt_2} d\mu(\phi)\EEQ as
\BEQ Z:=\int  e^{- \int_{\R} {\cal L}_{int}(\phi_1,\phi_2,\sigma)(t) dt}
d\mu(\phi) d\mu(\sigma),\EEQ where
 \BEQ {\cal L}_{int}(\phi_1,\phi_2,\sigma)(t)=\II\lambda\left(\partial A^{+}(t)\sigma_+(t)
-\partial {\cal A}^-(t)\sigma_-(t)\right). \label{eq:Lint} \EEQ

\medskip All of this is ill-defined mathematically since (1) $\sigma$ is a distribution-valued process and $\partial A^{\pm}$ is not defined {\em at all} when
$\alpha\le 1/4$; (2) one integrates over $\R$ a translation-invariant quantity (note
that $\phi_1,\phi_2,\sigma$ are all {\em stationary fields}).


\section{Heuristic perturbative proof of convergence}


Let us now explain the basics of perturbative quantum field theory, and show how it suggests (at least
heuristically) the assertions of Theorem 0.1.
\medskip
The general idea is  to expand formally the exponential of the Lagrangian in order to compute polynomial
moments, $\frac{1}{Z}
\esper\left[ \psi_1(x_1)\ldots\psi_n(x_n) e^{- \int {\cal L}_{int}(\phi_1,\phi_2,\sigma)(t) dt}
\right]$,
 also called {\em $n$-point functions} and
denoted by  $\langle \psi_1(x_1)\ldots\psi_n(x_n)\rangle_{\lambda}$, $\psi_i=\phi_1,\phi_2,
\sigma_+$ or $\sigma_-$, as $\frac{1}{Z}\sum_{n\ge 0} \frac{(-1)^n}{n!}
\esper\left[ \psi_1(x_1)\ldots\psi_n(x_n) \left( \int {\cal L}_{int}(\ .\ ;t) dt\right)^n \right]$. We do
not bother too much about the volume and ultra-violet cut-off here, and write $\langle \ \cdot \ \rangle_{\lambda}$
instead of $\langle \ \cdot \ \rangle_{\lambda,V,\rho}$.
Recall first the following classical combinatorial facts. A good reference for perturbative expansions in quantum
field theory is e.g. \cite{LeB}.

\begin{Proposition}
\begin{enumerate}
\item (Wick's formula) Let $X=(X_1,\ldots,X_{2n})$ be a (centered) Gaussian vector. Then
\BEQ \esper [X_1\ldots X_{2n}]=\sum_{(i_1 i_2)\ldots(i_{2n-1}i_{2n})} \esper[X_{i_1}X_{i_2}]\ldots\esper[X_{i_{2n-1}}X_{i_{2n}}],\EEQ
where the indices range over all  pairings of the indices $1,\ldots,2n$. Each term in the sum is represented as a graph with $2n$ points   connected  two by two.
\item(connected moments) Let $\langle \ \cdot\  \rangle:=\frac{\esper\left[ \ \cdot \ e^{\Phi(X)} \right]}{\esper[ e^{\Phi(X)}]}$ be a penalized measure, where $X=(X_1,X_2,\ldots)$ is a (centered) Gaussian
vector, and $\Phi(X)$ is a polynomial in $X_1,X_2,\ldots$. Then the connected expectation $\langle X_1\ldots X_n\rangle_c$ ($c$ for connected) is (formally at least) the sum of all
connected graphs obtained by (i) expanding the exponential; (ii) applying Wick's formula and drawing links between the paired points; (iii) identifying all points coming from the same
{\em vertex}, i.e. from the same monomial in $\Phi(X)$ descended from the exponential.
\end{enumerate}
\end{Proposition}

Connected expectations exclude in particular vacuum contributions, i.e. terms of the form $\esper[ e^{\Phi(X)}]
\esper [X_1\ldots X_n]=Z\esper[X_1\ldots X_n]$. Discarding these contributions
can be shown to provide automatically the normalizing factor $\frac{1}{Z}$. Then usual expectations $\langle X_1\ldots X_n\rangle$ are obtained by taking all possible splittings
of $\{1,\ldots,n\}$ into disjoint subsets $I_1\uplus\ldots\uplus I_p$ and summing over the products of connected expectations $\sum_p\sum_{I_1,\ldots,I_p}\langle \prod_{i\in I_1}
X_i\rangle_c\ldots \langle \prod_{i\in I_p} X_i\rangle_c$. In practice the last operation is trivial for two-point functions $\langle X_{i_1} X_{i_2}\rangle$ if by parity (which is often the case in quantum field
theory) the one-point functions $\langle X_i\rangle_c$ vanish, so that $\langle X_{i_1}X_{i_2}\rangle=
\langle X_{i_1}X_{i_2}\rangle_c$.

\bigskip

Let us return to our model.  Using a straightforward extension of the above Proposition, one may represent
 $\langle \psi_1(x_1)\ldots\psi_n(x_n)\rangle_{\lambda}$, $\psi=\phi$ or $\sigma$  as a sum over Feynman diagrams, $\sum_{\Gamma} A(\Gamma)$, where $\Gamma$ ranges over a set of diagrams with $n$
external legs, and $A(\Gamma)\in\R$ is the evaluation of the corresponding diagram (see examples below); connected expectations will then be obtained as a sum over connected Feynman
diagrams. More precisely, one obtains formally a (diverging) power series in $\lambda$, $\sum_{n\ge 0}\lambda^n \sum_{\Gamma_n} A(\Gamma_n)$, where $\Gamma_n$ ranges over the set of 
Feynman diagrams with $n$ vertices.   The Gaussian integration by parts formula \footnote{an infinite-dimensional
extension of the well-known formula for Gaussian vectors, $\esper \left[
\partial_{X_i}F(X_1,\ldots,X_n)\right]=\sum_j C^{-1}(i,j) \esper \left[ X_j F(X_1,\ldots,X_n)\right]$ if
$C$ is the covariance matrix of $(X_1,\ldots,X_n)$.} yields a so-called {\em Schwinger-Dyson identity},
\BEA  \langle \partial {\cal A}^{\pm}(x)\partial {\cal A}^{\pm}(y)\rangle_{\lambda}
&=&-\frac{1}{\lambda^2 Z(\lambda)} \esper\left[ \frac{\del}{\del \sigma_+(y)} \frac{\del}{\del\sigma_+(x)}
e^{-\int {\cal L}_{int}(\phi_1,\phi_2,\sigma_+)(t)dt} \right]\nonumber\\
&=& -\frac{1}{\lambda^2 Z(\lambda)} \esper\left[ (C_{\sigma_+}^{-1}\sigma_+)(y) \frac{\del}{\del\sigma_+(x)}
e^{- \int {\cal L}_{int}(\phi_1,\phi_2,\sigma_+)(t)dt} \right] \nonumber\\
&=& -\frac{1}{\lambda^2} \left[ -C_{\sigma_+}^{-1}(x,y)+ \langle (C_{\sigma_+}^{-1}\sigma_+)(x)
(C_{\sigma_+}^{-1}\sigma_+)(y)\rangle_{\lambda}\right], \nonumber\\ \EEA
with Fourier transform
\BEQ  \langle |{\cal F}(\partial {\cal A}^{\pm})(\xi)|^2\rangle_{\lambda}=\frac{1}{\lambda^2}
|\xi|^{1-4\alpha} \left[ 1-|\xi|^{1-4\alpha}  \langle |({\cal F}\sigma_+)(\xi)|^2\rangle_{\lambda} \right].
\label{eq:4.3}
\EEQ

By parity, $\langle |{\cal F}(\partial{\cal A}^{\pm})(\xi)|^2\rangle_{\lambda}$ is a power series in $\lambda^2$.

\begin{figure}[h]
  \centering
   \includegraphics[scale=0.3]{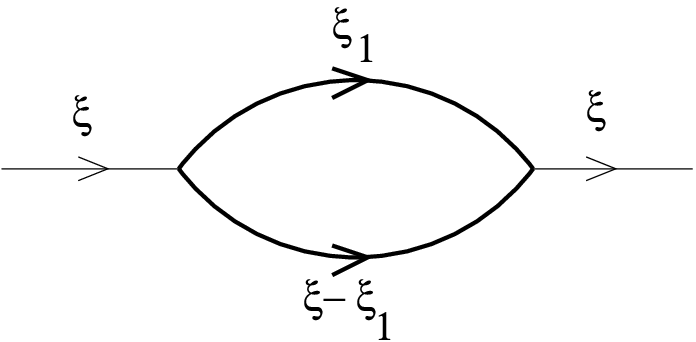}
   \caption{\small Bubble diagram with $2$ vertices. By momentum conservation $\xi=\xi_1+\xi_2$, which leaves out one free internal momentum.}
  \label{Fig-bubble}
\end{figure}

\begin{figure}[h]
  \centering
   \includegraphics[scale=0.3]{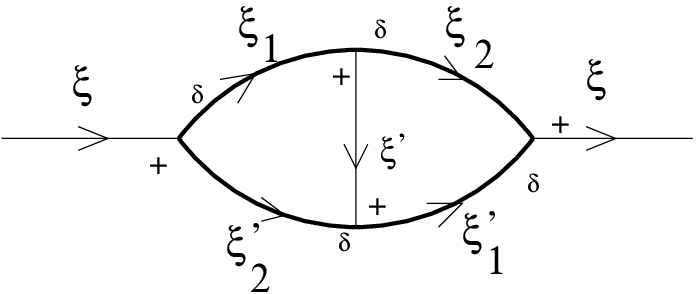}
   \caption{\small More complicated bubble diagram with $4$ vertices. By momentum conservation $\xi=\xi_1+\xi'_2=\xi'_1+\xi_2$ and $\xi_1=\xi'+\xi_2$, which leaves out two independent internal
   momenta.}
  \label{Fig-doublebubble}
\end{figure}

\begin{figure}[h]
  \centering
   \includegraphics[scale=0.4]{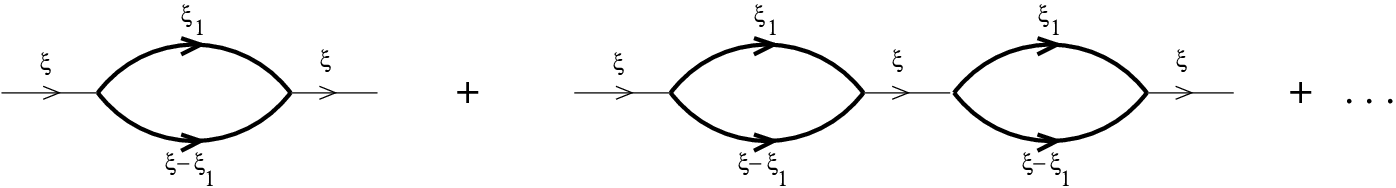}
   \caption{\small First three terms of the bubble series. The renormalized covariance of the $\sigma$-field
is equal to the sum of the series.}
  \label{Fig-bubble-series}
\end{figure}

Introduce an ultra-violet cut-off at scale $\rho$ as in Definition \ref{def:cut-off}. For the simplicity of the exposition we shall actually use a 
 brute-force ultraviolet cut-off at momentum $M^{\rho}$, i.e. cut off all Fourier components with momentum $|\xi|>M^{\rho}$. After Fourier transformation,
 $\int {\cal L}_{int}(\cdot;t)dt$ becomes $\II\lambda \int_{|\xi_1|<|\xi_2|} d\xi_1 d\xi_2 d\xi \del_0(\xi_1+\xi_2+\xi) {\cal F}\sigma_+(\xi) {\cal F}(\partial \phi_1)(\xi_1)
 {\cal F}\phi_2(\xi_2)$, minus a similar term involving $\sigma_-$. The square of this expression contributes the following term of order $O(\lambda^2)$ to $\langle |{\cal F}\sigma_+(\xi)|^2
 \rangle_{\lambda}$,
\BEA &&  (-\II\lambda)^2  \int^{M^{\rho}}_{|\xi_1|<|\xi-\xi_1|} d\xi_1
\left\{  \left( \esper[ |{\cal F}\sigma_+(\xi)|^2 ]  \right)^2 \
\esper[ |{\cal F}(\partial\phi_1)(\xi_1)|^2 ]\  \esper[ |{\cal F}\phi_2(\xi-\xi_1)|^2 ] \right\} \nonumber\\
&&=  -\lambda^2 |\xi|^{8\alpha-2}  \int_{|\xi_1|<|\xi-\xi_1|}^{M^{\rho}} d\xi_1 |\xi_1|^{1-2\alpha}
|\xi-\xi_1|^{-1-2\alpha} \sim_{\rho\to\infty} -K\lambda^2 |\xi|^{8\alpha-2}  (M^{\rho})^{1-4\alpha}. \nonumber\\  \label{eq:bubble} \EEA
This is the evaluation of the {\em Feynman diagram represented} in Fig. \ref{Fig-bubble}, according to the following rules.

\begin{Definition}[Feynman rules]

A Feynman diagram in our theory is made up of (1) {\em bold lines} of type $i=1,2$, with momenta $\xi_i,\xi'_i,\ldots$ evaluated as $\esper |{\cal F}\phi_i(\xi_i)|^2=\frac{1}{|\xi_i|^{1+2\alpha}}$;
(2) {\em plain lines} of type $\pm$, with momenta $\xi,\xi',\ldots$, evaluated as $\esper |{\cal F}\sigma_{\pm}(\xi)|^2=\frac{1}{|\xi|^{1-4\alpha}}$; (3) vertices where two plain lines -- one
of each type -- and a bold line meet, with a momentum conservation rule, $\xi=\pm\xi_1\pm\xi_2$ (depending on the orientation of the lines). The definition of the interaction implies
the presence of  a further derivation -- represented by the symbol $\partial$ on the Feynman diagram -- on the $\phi_1$-, resp. $\phi_2$-field, and a momentum scale restriction $|\xi_1|<|\xi_2|$, resp.
$|\xi_1|>|\xi_2|$, at vertices involving a $\sigma_+$-, resp. $\sigma_-$-field. The derivation translates into a multiplication by $\II \xi_1$, resp. $\II \xi_2$ when evaluating the diagram.

\end{Definition}

It is sometimes useful to consider the evaluation of the corresponding {\em amputated Feynman diagram}, from
which  the contribution of the external legs has been removed.  Here for instance, the
evaluation of the amputated Feynman diagram associated to Fig. \ref{Fig-bubble} is $(|\xi|^{1-4\alpha})^2$ times the previous expressions, hence is equivalent to
the $\xi$-independent expression $-K\lambda^2 (M^{\rho})^{1-4\alpha}$ when $\rho\to\infty$.
It is a diverging {\em negative} quantity. (Using the Fourier truncation of Definition
\ref{def:cut-off} only changes the constant $K$.)  However, resumming {\em formally} the bubble series as in
Fig. \ref{Fig-bubble-series}
yields, starting from the right-hand side of eq. (\ref{eq:4.3}),  
\BEA &&  \frac{1}{\lambda^2}|\xi|^{1-4\alpha} \left[1-\sum_{n\ge 0} (-1)^n \left(\frac{1}{|\xi|^{1-4\alpha}}\cdot K\lambda^2 (M^{\rho})^{1-4\alpha}) \right)\right] \nonumber\\
&& \qquad =\frac{1}{\lambda^2}|\xi|^{1-4\alpha}  \ \cdot\
\frac{K\lambda^2 (M^{\rho}/|\xi|)^{1-4\alpha}}{1+K\lambda^2 (M^{\rho}/|\xi|)^{1-4\alpha}} \nonumber\\
&&\qquad \to_{\rho\to\infty}  \frac{1}{\lambda^2} |\xi|^{1-4\alpha}. \label{eq:bubble-series} \EEA
 
On the other hand (see Fig. \ref{Fig-bubble-series}),  the bare $\sigma$-covariance
 $\frac{1}{|\xi|^{1-4\alpha}}$ has been replaced with the renormalized covariance
 \BEQ\frac{1}{|\xi|^{1-4\alpha}} \cdot \frac{1}{1+K\lambda^2 (M^{\rho}/|\xi|)^{1-4\alpha}}=
\frac{1}{|\xi|^{1-4\alpha}+K\lambda^2 (M^{\rho})^{1-4\alpha}}, \EEQ which vanishes in the
 limit $\rho\to\infty$. The essential reason for this is of course that the oscillating signs
$(-1)^n$ in the bubble series evaluation -- due to the fact that the interaction Lagrangian
${\cal L}_{int}(\phi_1,\phi_2,\sigma)$ is purely {\em imaginary} -- result by summing in a huge,
virtually infinite denominator. Taking into account the possible insertion of $\sigma_-$-lines between $\sigma_+$-lines
 amounts to a simple change of the constant $K$.
In physical terms, the interaction in $\frac{1}{|\xi|^{1-4\alpha}}$ has been  {\em screened} by a {\em huge
mass term} $K\lambda^2 M^{\rho(1-4\alpha)}\to_{\rho\to\infty} +\infty$ (see section 3 for the definition of the mass). More complicated diagrams contributing to $\langle |({\cal F}\sigma_+)(\xi)|^2\rangle_{\lambda}$, and involving internal $\sigma$-lines as in Fig. \ref{Fig-doublebubble} also vanish when $\rho\to\infty$. Thus there remains simply:

\BEQ \langle |{\cal F} {\cal A}^{\pm}(\xi)|^2|\rangle_{\lambda}=\frac{1}{\lambda^2}
|\xi|^{-1-4\alpha}. \EEQ
Hence $\esper |{\cal A}^{\pm}(t)-{\cal A}^{\pm}(s)|^2 \lesssim \frac{1}{\lambda^2} |t-s|^{4\alpha}$, as in
eq. (\ref{eq:1.13}).

\medskip

As for the mixed term $\langle \partial {\cal A}^{\pm}(x)\partial{\cal A}^{\mp}(y)\rangle_{\lambda}$, its
Fourier transform is given by $\frac{1}{\lambda^2}|\xi|^{1-4\alpha} \left[-\frac{1}{1+K''\lambda^2 (\Lambda/|\xi|)^{1-4\alpha}}\right]$, where $K''<K$ due to the constraints on the scales for
bubbles of mixed type with one $\sigma_+$- and one $\sigma_-$-leg, which vanishes
in the limit $\rho\to\infty$ (note the disappearance of the factor $1$ compared to eq. (\ref{eq:bubble-series}), due to the fact that $\esper \sigma^+(x)\sigma^-(y)=0$). Thus the covariance of the two-component
$\sigma$-field has been renormalized to $\frac{1}{|\xi|^{1-4\alpha}\Id+m^{\rho}}$, where $m^{\rho}$ is
a two-by-two positive "mass" matrix with eigenvalues $\thickapprox \lambda^2 M^{\rho(1-4\alpha)}$.

\medskip

Using eq. (\ref{eq:1.8}), one obtains:
\BEA (2\pi c_{\alpha})^2 \langle {\cal A}(s,t)^2 \rangle_{\lambda} &=& \langle \left| {\cal A}^+(t)-
{\cal A}^+(s)\right|^2 \rangle_{\lambda}+\langle \left| {\cal A}^-(t)-
{\cal A}^-(s)\right|^2 \rangle_{\lambda} \nonumber\\
&& \qquad \qquad \qquad +\esper \left|{\cal A}^+_{\partial}(s,t)-
{\cal A}_{\partial}^-(s,t)\right|^2 \nonumber\\
&=&   \frac{4}{\lambda^2} \int (1-\cos(t-s)\xi)|\xi|^{-1-4\alpha} d\xi+
 \esper \left|{\cal A}^+_{\partial}(s,t)-
{\cal A}_{\partial}^-(s,t)\right|^2 \nonumber\\
&=& (\frac{4}{\lambda^2}K_1+K_2)|t-s|^{4\alpha} \EEA
for some constants $K_1,K_2$.

\bigskip

Let us now consider briefly other correlations. For a general discussion we need the following easy power-counting lemma:

\begin{Lemma}[power-counting rules]

Let $\Gamma$ be a Feyman diagram with $N_{\sigma}$ external $\sigma$-lines, $N_{\phi}$ external $\phi$-lines, and $N_{\partial\phi}$ external $\partial \phi$-lines. Then the
overall degree of homogeneity (in powers of $\xi$) of the evaluation of the corresponding amputated diagram 
-- also called: {\em overall degree of divergence} --  is ${1-2\alpha N_{\sigma}+\alpha N_{\phi}+(\alpha-1)N_{\partial\phi}}$.

\end{Lemma}

{\bf Proof.} Let: $I_{\sigma}$, resp. $I_{\phi}$, be the number of {\em internal} lines of type $\sigma$,
resp.  $\phi$ or  $\partial\phi$; $I=I_{\sigma}+I_{\phi}$ be the total number of internal lines;
and $L=I-V+1$ be the number of loops, equal to the number of independent momenta (one per internal line, minus one per vertex due to momentum conservation, plus one due to overall momentum
conservation). Since one $\sigma$- and two $\phi$-lines meet at each vertex, one also has the relations $2I_{\sigma}+N_{\sigma}=V$, and $2I_{\phi}+N_{\phi}+N_{\partial\phi}=2V$. Now  the amputated diagram is
homogeneous to  $|\xi|^{-(1-4\alpha)I_{\sigma}-(1+2\alpha)I_{\phi}+L+V-N_{\partial\phi}}$ (counting one derivative per vertex, and minus one derivative
per {\em external} $\partial\phi$-leg which is not taken into account in the evaluation). Putting all these relations together yields the result. \hfill \eop

If a diagram is overall divergent, i.e. if its overall degree of divergence  is positive, then the diagram diverges (except if by chance the coefficient of the term of highest degree in $\xi$ vanishes).
 On the other hand, the fact that a diagram is overall convergent (i.e. its overall degree of divergence is negative) does not imply that it is convergent, since it may contain overall divergent
 {\em sub}-diagrams. One must hence study the behaviour of all possible  diagrams, with arbitrary external leg structure.

\begin{figure}[h]
  \centering
   \includegraphics[scale=0.4]{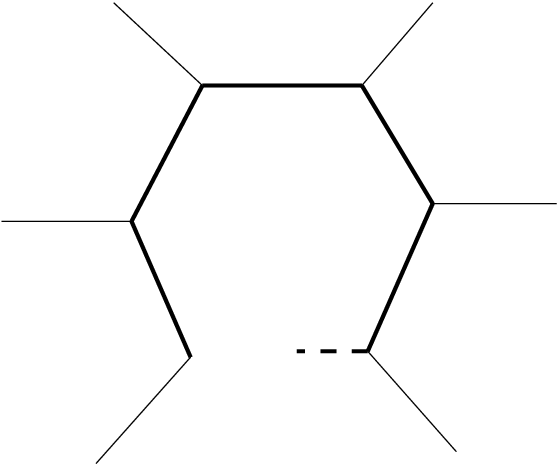}
   \caption{Higher connected moments of the L\'evy area.}
  \label{Fig-polygonal}
\end{figure}
 
The above simple power-counting argument shows that the overall degree of divergence of a connected diagram with $2n$ external $\sigma$-legs is
$1-4n\alpha$. For $n\ge 2$, this is $\le 1-8\alpha<0$ since $\alpha>\frac{1}{8}$ by hypothesis, so such diagrams are overall convergent. By the above arguments, there remain only the connected diagrams in the limit $\Lambda\to\infty$,  see Fig. \ref{Fig-polygonal}, whose evaluation is independent of $\lambda$.

General considerations following from the multi-scale expansions (one may refer to \cite{FVT} for a good, accessible presentation, or to \cite{Unt-ren} for an application to the Gaussian renormalization
of iterated integrals evoked in subsection 2.3) show that it is enough to consider the behaviour of diagrams whose internal legs have higher (or even: much higher) momentum scale than external legs, the so-called
{\em dangerous diagrams}. Then the
momentum scale constraint on the vertices coming from Fourier normal ordering implies that  the external legs of dangerous diagrams may be either of type $\sigma$ or of type $\partial\phi$,
but not of type $\phi$. Consider now any diagram whose external structure contains external $\partial\phi$-legs. By parity it has at least two such external legs, and the previous power-counting
rules show that such a diagram is always overall convergent.

Finally,  the law of the field $\phi$ is left unchanged by the interaction. Namely,
 all non-trivial diagrams contributing e.g.  to $\langle \phi_1(x)\phi_2(x)\rangle_{\lambda}$ involve internal $\sigma$-lines which (as previously "shown") vanish in the limit $\rho\to\infty$.

On the whole, this is the content of Theorem 0.1.

\bigskip

The art of constructive field theory is to make the previous speculations rigorous. It relies on the following considerations, corresponding to the weak points (not to say flaws!) in the
above arguments:

\begin{enumerate}
\item While going from eq. (\ref{eq:bubble}) to (\ref{eq:bubble-series}), we have replaced the amputated bubble diagram evaluation by its asymptotics when $\rho\to\infty$, namely,
$-K\lambda^2 M^{\rho(1-4\alpha)}$, which is simply equal to its {\em evaluation at zero external momentum $\xi$}, also called {\em local part}. Thus we have actually not resummed the
whole bubble series, but only the corresponding local parts, and observed that this was equivalent to adding a  mass term of the form  $K'\lambda^2 M^{\rho(1-4\alpha)} \int |\sigma(x)|^2 dx$ to
the Lagrangian.

\item The bubble series is really a terribly {\em diverging geometric series}. Renormalization must actually be performed scale by scale. Considering only bubble diagrams with momentum in the
dyadic slice $M^{\rho-1}<|\xi|<M^{\rho}$ leads on the other hand to a converging geometric series for $\lambda$ small enough since the term between parentheses in eq. (\ref{eq:bubble-series}), $K\lambda^2 \left( \frac{M^{\rho}}{|\xi|}\right)^{1-4\alpha}$,  is then $<1$.  This is equivalent to integrating out
the highest field components $(\sigma^{\rho},\phi^{\rho})$, as explained in section 3. One obtains thus a running mass coefficient $m^{\rho}$ of order $\lambda^2 M^{(1-4\alpha)\rho}$.
The procedure must then be iterated by going down the scales step by step. Since  renormalization
{\em reduces} the  covariance of the $\sigma$-field, the
bound on $\lambda$ ensuring convergence does not become worse and worse after each step.

\item We neglected more complicated bubble diagrams as in Fig. \ref{Fig-doublebubble}. Although these have the same order as the simple bubble diagram of Fig. \ref{Fig-bubble}, as follows
from the above power-counting rules, taking into consideration {\em all possible} bubble diagrams lead to a terribly diverging power series in $\lambda$ due to the rapidly increasing
number of such diagrams in terms of the number of vertices, with a coefficient roughly of order $n!$ in front of $\lambda^n$. This divergence is actually due to the accumulation of vertices
in a small region of space of size $O(M^{-j})$, where $j$ is the momentum scale under consideration. {\em Multi-scale cluster  expansions} in constructive field theory, by considering only {\em partial} series expansions,
avoid this dangerous accumulation process.

\item By splitting each vertex $\int {\cal L}_{int}^{\to\rho}(\cdot;x)dx$ into its different scales,
 there may appear fields  $\phi_1^{j_1},\phi_2^{j_2},\sigma^j$ with different scales $j_1\not=j_2\not=j$. Taking
this into account in a coherent way in the previous partial series expansions lead to complicated combinatorial expressions encoded by so-called {\em polymers}, which are the main
object in use in constructive field theory.

\item In the previous vertex splitting, the field with lowest momentum scale ($j_1$, $j_2$ or $j$, depending on the case) is called {\em low-momentum field}. Even though the cluster expansion
in each momentum scale prevents an accumulation of vertices in the same region of space, the compound effect of {\em all} cluster expansions at {\em all} scales produces unavoidably
accumulations of fields with very low momentum in very large regions of space, which is a dangerous problem called {\em domination problem}. This accounts for the addition of the 
extra boundary term ${\cal L}_{bdry}^{\to\rho}$ in the interaction Lagrangian. Writing out this term and explaining its precise form would however take us too far away.

\end{enumerate}


\end{document}